\documentclass[11pt]{amsart} \usepackage{amscd, amssymb}

\begin{document}
\def\LU{{\mathcal L\mathcal U}} \def\dach{\!\widehat{\phantom{G}}}
\def\Aut{\operatorname{Aut}} \def\Prim{\operatorname{Prim}}
\def\Inn{\operatorname{Inn}} \def\Ad{\operatorname{Ad}}
\def\id{\operatorname{id}} \def\supp{\operatorname{supp}}

\def\Spec{\operatorname{Spec}} \def\Ind{\operatorname{Ind}}
\def\nd{\mathrm{es}} \def\triv{\mathrm{triv}} \def\eps{\varepsilon}
\def\K{{\mathcal K}} \def\L{\mathcal L} \def\R{\mathcal R}
\def\L{\mathcal L} \def\C{\mathcal C} \def\E{\mathcal E}
\def\Q{\mathcal Q} \def\F{\mathcal F} \def\B{\mathcal B}
\def\H{\mathcal H} \def\V{\mathcal V} \def\lk{\langle}
\def\rk{\rangle} \def\EE{\mathbb E} \def\DD{\mathbb D}
\def\I{{\operatorname{I}}} \def\U{{\mathcal U}} \def\UM{{\mathcal U}M}
\def\ZUM{{\mathcal Z}UM} \def\ZM{{\mathcal Z}M}
\def\KK{\operatorname{KK}} \def\RK{\operatorname{RK}}
\def\RKK{\operatorname{RKK}} \def\ind{\operatorname{ind}}
\def\res{\operatorname{res}} \def\inf{\operatorname{inf}}
\def\ker{\operatorname{ker}} \def\infl{\operatorname{inf}}
\def\pt{\operatorname{pt}} \def\M{\mathfrak M}
\def\comp{\operatorname{comp}}\def\infl{\operatorname{infl}}
\def\Bott{\operatorname{Bott}} \def\GL{\operatorname{GL}}
\def\sp{\operatorname{span}} \def\k{\operatorname{K}}
\def\EEG{\underline{\underline{\mathcal{E}G}}} \def\EG{\mathcal{E}(G)}
\def\EH{\mathcal{E}(H)} \def\EGT{\mathcal{E}(\tilde{G})}
\def\EGN{\mathcal{E}(G/N)} \renewcommand{\top}{\operatorname{top}}
\def\A{\mathcal{A}} \def\sm{\backslash} \def\Ind{\operatorname{Ind}}
\def\TT{\mathbb T} \def\ZZ{\mathbb Z} \def\CC{\mathbb C}
\def\FF{\mathbb F} \def\RR{\mathbb R} \def\NN{\mathbb N}
\def\om{\omega} \def\Spin{\operatorname{Spin}}
\def\Pin{\operatorname{Pin}}

\newcommand{\ddS}{\stackrel{\scriptscriptstyle{o}}{S}}
\renewcommand{\oplus}{\bigoplus} \newtheorem{thm}{Theorem}[section]
\newtheorem{cor}[thm]{Corollary} \newtheorem{ass}{Assumption}[section]
\newtheorem{lem}[thm]{Lemma} \newtheorem{prop}[thm]{Proposition}
\newtheorem{lemdef}[thm]{Lemma and Definition}

\theoremstyle{definition} \newtheorem{defn}[thm]{Definition}

\theoremstyle{remark} \newtheorem{remark}[thm]{\bf Remark}
\newtheorem{ex}[thm]{\bf Example}

\numberwithin{equation}{section} \emergencystretch 30pt

\renewcommand{\theenumi}{\roman{enumi}}
\renewcommand{\labelenumi}{(\theenumi)}

\title[Connes-Kasparov Conjecture]{The Connes-Kasparov Conjecture for
  almost connected groups}

\author[Chabert]{J\'er\^ome Chabert} \author[Echterhoff]{Siegfried
  Echterhoff} \author[Nest]{Ryszard Nest} \address{J\'er\^ome Chabert:
  Universit\'e Blaise Pascal, B\^at. de math\'ematiques, 63177,
  Aubi\`ere, France} \address{Westf\"alische Wilhelms-Universit\"at
  M\"unster, Mathematisches Institut, Einsteinstr. 62 D-48149
  M\"unster, Germany} \address{Department of Mathematics, University
  of Copenhagen, Universitetsparken 5, DK-2100 Copenhagen, Denmark}
\email{chabert@math.univ-bpclermont.fr, echters@math.uni-muenster.de,
  rnest@math.ku.dk} \thanks{This research has been supported by the
  Deutsche Forschungsgemeinschaft (SFB 478)}
\begin{abstract}
  Let $G$ be a locally compact group with cocompact connected
  component. We prove that the assembly map from the topological
  $\k$-theory of $G$ to the $\k$-theory of the reduced $C^*$-algebra
  of $G$ is an isomorphism.
\end{abstract}

\maketitle
\tableofcontents

\section{Introduction and statement of results}\label{sec-intro}

In this paper we give a proof of the Connes-Kasparov conjecture for
almost connected groups.  To be more precise, we prove the following

\begin{thm}\label{thm-main}
  Let $G$ be a second countable almost connected group (i.e., $G/G_0$
  is compact, where $G_0$ denotes the connected component of $G$).
  Then $G$ satisfies the Baum-Connes conjecture with trivial
  coefficients $\CC$, i.e., if $ \k_*^{\top}(G)$ denotes the
  topological $K$-theory of $G$, then the Baum-Connes assembly map
  $$\mu: \k_*^{\top}(G)\to \k_*(C_r^*(G))$$
  is an isomorphism.
\end{thm}

It was already shown by Kasparov in \cite{Kas1} that the theorem is
true if $G$ is amenable. In fact, by a more recent result of Higson
and Kasparov, we know that the Baum-Connes conjecture with arbitrary
coefficients holds for any amenable group.  By work of A. Wassermann
\cite{Was}, we also know that the result is true for all connected
reductive linear Lie groups. More recently, Lafforgue used quite
different methods to give a proof of the conjecture for all connected
semi-simple groups with finite center (which are not necessarily
linear).  The main idea of the proof of the general result of Theorem
\ref{thm-main} is to use the Mackey-machine approach, as outlined in
\cite{CE2}, in order to reduce to the reductive case.  The strategy
for doing this bases heavily on some ideas presented in Puk\'anszky's
recent book \cite{Puk} where he reports on his deep analysis of the
representation theory of connected groups. In particular the methods
of his proof that locally algebraic connected real Lie groups are type
I, presented on the first four pages of his book, where most
enlightening.

The above theorem is actually a special case of a more general result
which we shall explain below. If $G$ is a second countable locally
compact group, then by a $G$-algebra $A$ we shall always understand a
$C^*$-algebra $A$ equipped with a strongly continuous action of $G$ by
$*$-automorphisms of $A$. Let $\E(G)$ denote a locally compact
universal proper $G$-space in the sense of \cite{KS} (we refer to
\cite{CEM} for a discussion about the relation to the notion of
universal proper $G$-space as introduced by Baum, Connes and Higson in
\cite{BCH}).  If $A$ is a $G$-algebra, the topological $K$-theory of
$G$ with coefficients in $A$ is defined as
$$\k_*^{\top}(G,A)=\lim_{X}\KK_*^G(C_0(X),A),$$
where $X$ runs through
the $G$-compact subspaces of $\E(G)$ (i.e., $X/G$ is compact) ordered
by inclusion, and $\KK_*^G(C_0(X),A)$ denotes Kasparov's equivariant
$\KK$-theory. If $A=\CC$, we simply write $\k_*^{\top}(G)$ for
$\k_*^{\top}(G,\CC)$.

The construction of Baum, Connes and Higson presented in \cite[\S
9]{BCH} determines a homomorphism
$$\mu_A:\k_*^{\top}(G,A)\to \k_*(A\rtimes_rG),$$
usually called the
{\em assembly map}. We say that $G$ satisfies BC for $A$ (i.e., $G$
satisfies the Baum-Connes conjecture for the coefficient algebra $A$),
if $\mu_A$ is an isomorphism.  Theorem \ref{thm-main} is then a
special case of

\begin{thm}\label{thm-main1}
  Suppose that $G$ is any second countable locally compact group such
  that $G/G_0$ satisfies BC for arbitrary coefficients, where $G_0$
  denotes the connected component of $G$.  (By the results of Higson
  and Kasparov \cite{HK} this is in particular true if $G/G_0$ is
  amenable or, more general, if $G/G_0$ satisfies the Haagerup
  property.)  Then $G$ satisfies BC for $\K(H)$, $H$ a separable
  Hilbert space, with respect to any action of $G$ on $\K(H)$.
\end{thm}

It is well known that in case of almost connected groups, the
topological $K$-theory $\k_*^{\top}(G,A)$ has a very nice description
in terms of the maximal compact subgroup $L$ of $G$. In fact, under
some mild extra conditions on $G$, the group $\k_*^{\top}(G,A)$ can be
computed by means of the $K$-theory of the crossed product $A\rtimes
L$. We give a brief discussion of these relations in \S\ref{sec-max}
below. As was already pointed out in \cite{ros}, our results have
important applications to the study of square-integrable
representations.  In fact, combining our results with \cite[Theorem
4.6]{ros} gives

\begin{cor}[{cf \cite[Corollary 4.7]{ros}}]\label{cor-square0}
  Let $G$ be a connected unimodular Lie group.  Then all
  square-integrable factor representations of $G$ are type I.
  Moreover, $G$ has no square-integrable factor representations if
  $\dim(G/L)$ is odd, where $L$ denotes the maximal compact subgroup
  of $G$.
\end{cor}

The paper is outlined as follows: In our preliminary section,
\S\ref{sec-prel}, we recall the main results from \cite{CE1,CE2} on
the permanence properties of the Baum-Connes conjecture which are
needed in this work. We will also use these results to perform some
first reductions of the problem. In \S\ref{sec-fields} we prove a
result on continuous fields of actions, showing under some mild
conditions on the group $G$ and the base space $X$ of the field, that
$G$ satisfies the Baum-Connes conjecture with coefficients in the
algebra of $C_0$-sections of the field if it satisfies the conjecture
for all fibres. This result will be another basic tool for the proof
our our main theorem.

In \S\ref{sec-semi} we are concerned with the conjecture for reductive
groups. Using (and slightly extending) some recent results of
Lafforgue \cite{Laf} on the conjecture for semi-simple groups with
finite center, we will show that the results on continuous fields
obtained in \S\ref{sec-fields} imply that the conjecture (with trivial
coefficients) holds for all reductive groups without any linearity
conditions.  In \S\ref{sec-perfect} we will then use Puk\'anzsky's
methods in combination with an extensive use of the permanence
properties for BC to give the final steps for the proof of Theorem
\ref{thm-main1}.  In \S\ref{sec-max} we shall then discuss the
implications for almost connected groups as mentioned above. Note that
\S\ref{sec-max} does not contain any new material, except the
conclusion drawn out of our main theorem.

{\bf Acknowledgments.} The authors profited quite a lot from valuable
discussions with several friends and colleagues.  We are most grateful
to Alain Connes who suggested (after an Oberwolfach lecture on the
permanence properties for BC) to use the results of \cite{CE1, CE2} to
attack the Connes-Kasparov conjecture with trivial coefficients.  
We are also grateful to Bachir Bekka, Ludwig Br\"ocker, Guennadi
Kasparov, Guido Kings, Herv\'e Oyono-Oyono, J\"org Sch\"urmann 
and Peter Slodowy for some useful
discussions or comments.

\section{Some preliminaries and first reductions}\label{sec-prel}

Let us collect some general facts which were presented in \cite{CE2}
-- for the definitions of twisted actions and twisted equivariant
$\KK$-theory we refer to \cite{CE1}.  Assume that $G$ is a second
countable group and let $B$ be a $G$-$C^*$-algebra. We say that $G$
satisfies BC with coefficients in $B$ if the assembly map
$$\mu_B: \k_*^{\top}(G,B)\to \k_*(B\rtimes_rG)$$
is an isomorphism.
If $N$ is a closed normal subgroup of $G$, there exists a twisted
action of $(G,N)$ on $B\rtimes_rN$ such that the twisted crossed
product $(B\rtimes_rN)\rtimes_r(G,N)$ is canonically isomorphic to
$B\rtimes_rG$. Moreover, we can use the twisted equivariant
$\KK$-theory of \cite{CE1} to define the topological $\k$-theory $
\k_*^{\top}(G/N, B\rtimes_rN)$ with respect to the twisted action of
$(G,N)$ on $B\rtimes_rN$, and a twisted version of the assembly map
$$\mu_{B\rtimes_rN}: \k_*^{\top}(G/N, B\rtimes_rN)\to
\k_*((B\rtimes_rN)\rtimes_r(G,N)).$$
In \cite{CE1} we constructed a
partial assembly map
$$\mu_{N,B}^G: \k_*^{\top}(G,B)\to \k_*^{\top}(G/N, B\rtimes_rN)$$
such that the following diagram commutes
$$
\begin{CD}
  \k_*^{\top}(G,B) @>\mu_{N,B}^G>>  \k_*^{\top}(G/N, B\rtimes_rN)\\
  @V\mu_BVV   @VV\mu_{B\rtimes_rN}V\\
  \k_*(B\rtimes_rG) @>\cong>> \k_*((B\rtimes_rN)\rtimes_r(G,N)).
\end{CD}
$$
Using this, the first two authors were able to prove the following
extension results:

\begin{thm}\label{thm-ext}
  Assume that $B$ is a $G$-algebra and let $N$ be a closed normal
  subgroup of $G$. Let $q:G\to G/N$ denote the quotient map and assume
  that one of the following conditions is satisfied
      \begin{enumerate}
      \item $G/N$ has a compact open subgroup $\dot{K}$ and for any
        compact subgroup $\dot{C}$ of $G/N$, the group
        $C=q^{-1}(\dot{C})$ satisfies BC for $B$.
      \item $G$ has a $\gamma$-element $\gamma\in \KK_0^G(\CC,\CC)$
        (which is automatically true if $G$ is almost connected),
        $G/N$ is almost connected and $K=q^{-1}(\dot{K})$ satisfies BC
        for $B$, where $\dot{K}$ is a maximal compact subgroup of
        $G/N$.
      \end{enumerate}
      Then the partial assembly map $\mu_{N,B}^G: \k_*^{\top}(G,B)\to
      \k_*^{\top}(G/N, B\rtimes_rN)$ is an isomorphism. In particular,
      $G$ satisfies BC for $B$ if and only if $G/N$ satisfies BC for
      $B\rtimes_rN$.
\end{thm}
\begin{proof} See \cite[Theorem 3.3 and Theorem 3.7]{CE2}.
\end{proof}

In order to avoid the use of twisted actions we may use the version of
the Packer-Raeburn stabilization trick as given in \cite{PR, E-mor}:

\begin{prop}[{cf \cite[Theorem 3.4]{PR} and \cite[Corollary
    1]{E-mor}}]\label{prop-mor}
  Assume that $G$ is a second countable group and let $N$ be a closed
  normal subgroup of $G$.  Let $(\alpha,\tau)$ be a twisted action of
  $(G,N)$ on the separable $C^*$-algebra $A$. Then there exists an
  ordinary action $\beta:G/N\to \Aut(A\otimes \K)$, $\K=\K(l^2(\NN))$,
  such that $\beta$ is stably exterior equivalent (and hence Morita
  equivalent) to $(\alpha,\tau)$.
\end{prop}

Note that BC is invariant under passing to Morita equivalent actions.
Thus, in order to conclude that $(G,N)$ satisfies BC for
$B\rtimes_rN$, it is enough to show that $G/N$ satisfies BC for
$(B\rtimes_rN)\otimes\K$ with respect to an appropriate action of
$G/N$ on $(B\rtimes_rN)\otimes\K$. In particular, if $G/N$ is
amenable, it follows that $\mu_{B\rtimes_rN}: \k_*^{\top}(G/N,
B\rtimes_rN)\to \k_*((B\rtimes_rN)\rtimes_r(G,N))$ is always an
isomorphism.

In what follows we need to study the following special situation:
Assume that $\alpha:G\to \Aut(\K)$ is an action of $G$ on $\K=\K(H)$
for some separable Hilbert space $H$. Since $\Aut(\K)\cong
PU(H)=U(H)/\TT1$, we can choose a Borel map $V:G\to U(H)$ such that
$\alpha_s=\Ad V_s$ for all $s\in G$.  Since $\alpha$ is a
homomorphism, we see that there exists a Borel cocycle $\om\in
Z^2(G,\TT)$ such that
$$V_sV_t=\om(s,t)V_{st}\quad\text{for all $s,t\in G$}.$$
The class
$[\om]\in H^2(G,\TT)$ is called the {\em Mackey obstruction} for
$\alpha$ being unitary.  Let
$$1\mapsto \TT\to G_{\om}\to G\to 1$$
be the central extension of $G$
by $\TT$ corresponding to $\om$, i.e., we have $G_{\om}=G\times \TT$
with multiplication given by
$$(g,z)(g', z')=(gg', \om(g,g')zz'),$$
and the unique locally compact
group topology which generates the product Borel structure on $G\times
\TT$ (see \cite{Mac}).  Then the following is true

\begin{lem}\label{lem-twist}
  For each $n\in \ZZ$ let $\chi_n:\TT\to \TT; \chi_n(z)=z^n$.  Let
  $\alpha:G\to \Aut(\K)$ and $G_{\om}$ be as above.  Then $\alpha$ is
  Morita equivalent to the twisted action $(\id, \chi_1)$ of
  $(G_{\om},\TT)$ on $\CC$.
\end{lem}
\begin{proof}
  Let $V:G\to U(H)$ be as in the discussion above, i.e., $\alpha_s=\Ad
  V_s$ and $V_sV_t=\om(s,t)V_{st}$ for all $s,t\in G$. Then it is easy
  to check that $\tilde{V}:G_{\om}\to U(H)$ defined by
  $\tilde{V}_{(s,z)}=zV_s$ is a homomorphism which implements the
  desired equivalence on the $\K-\CC$ bimodule $H$ (we refer to
  \cite{E-mor} for an extensive discussion of Morita equivalence for
  twisted actions).
\end{proof}

Another important result is the continuity of the Baum-Connes
conjecture with respect to inductive limits of the coefficients, at
least if $G$ is exact. For this we need

\begin{lem}\label{lem-exactlimit}
  Assume that $(B_i)_{i\in I}$ is an inductive system of $G$-algebras
  and let $B=\lim_iB_i$ be the $C^*$-algebraic inductive limit.
  Assume further that one of the following conditions is satisfied:
\begin{enumerate}
\item All connecting maps $B_i\to B_j$, $i\leq j\in I$ are injective,
  or
\item $G$ is exact.
\end{enumerate}
Then $B\rtimes_rG=\lim_i(B_i\rtimes_rG)$ with respect to the obvious
connecting homomorphisms.
\end{lem}
\begin{proof} If all connecting maps are injective, we may
  regard each $B_i$ as a subalgebra of $B$.  But this implies that we
  also have $B_i\rtimes_rG$ as subalgebras of $B\rtimes_rG$, and hence
  the inductive limit
  $\lim_i(B_i\rtimes_rG)=\overline{\cup\{B_i\rtimes_rG: i\in I\}}$
  sits inside $B\rtimes_rG$. But it is easy to check that
  $\cup\{C_c(G,B_i):i\in I\}\subseteq \lim(B_i\rtimes_rG)$ is dense in
  $B\rtimes_rG$.
  
  Suppose now that $G$ is exact. In this situation we want to reduce
  the proof to situation (i).  Consider the canonical homomorphisms
  $\phi_i:B_i\to B$.  Let $I_i=\ker\phi_i$ and let
  $I_{ij}=\ker\phi_{ij}$, where the $\phi_{ij}:B_i\to B_j$ denote the
  connecting homomorphisms for $j\geq i$.  Of course, if $i\leq j\leq
  j'$ then $I_{ij}\subseteq I_{ij'}$, so for each $i\in I$ the system
  $(I_{ij})_{j\geq i}$ is an inductive system with injective
  connecting maps. It follows directly from the definition of the
  inductive limit that $I_i=\overline{\cup\{I_{ij}:j\geq
    i\}}=\lim_{j\geq i}I_{ij}$, and hence it follows from (i) that
  $I_i\rtimes_rG=\lim_{j\geq i}(I_{ij}\rtimes_rG)$.  By exactness of
  $G$ it follows that $I_i\rtimes_rG$ is the kernel of
  $\phi_i\rtimes_rG:B_i\rtimes_rG \to B\rtimes_rG$.  By the previous
  discussion it follows that $I_i\rtimes_rG=\lim_{j\geq
    i}(I_{ij}\rtimes_rG)$ is also the kernel of the canonical
  homomorphism $B_i\rtimes_rG\to \lim(B_j\rtimes_rG)$.  Thus, dividing
  out the kernels, i.e., by considering the system $(B'_i)_{i\in I}$
  with $B_i'=B_i/I_i$ we conclude from another use of (i) that
  $$B\rtimes_rG=\lim(B'_i\rtimes_rG)=\lim(B_i\rtimes_rG).$$
\end{proof}

As a direct consequence we obtain

\begin{prop}\label{prop-limit}
  Assume that the $G$-algebra $B$ is an inductive limit of the
  $G$-algebras $B_i$, $i\in I$, such that $G$ satisfies BC for all
  $B_i$.  Assume further that $G$ is exact or that all connecting
  homomorphisms $B_i\to B_j$ are injective.  Then $G$ satisfies BC for
  $B$.
\end{prop}
\begin{proof}
  It follows from Lemma \ref{lem-exactlimit} and the continuity of
  $K$-theory that $ \k_*(B\rtimes_rG)=\lim_i \k_*(B_i\rtimes_rG)$.  On
  the other side, it is shown in \cite[Proposition 7.1]{CE2} that
  $\k_*^{\top}(G,B)=\lim_i \k_*^{\top}(G, B_i)$.  Since by assumption
  $ \k_*^{\top}(G,B_i)\cong \k_*(B_i\rtimes_rG)$ via the assembly map,
  and since the assembly map commutes with the $K$-theory maps induced
  by the $G$-equivariant homomorphism $B_i\to B_j$, the result
  follows.
\end{proof}

As a first application we get

\begin{prop}\label{prop-central}
  Let $G$ be a separable locally compact group such that $G/G_0$
  satisfies BC for arbitrary coefficients.  Then the following are
  equivalent:
      \begin{enumerate}
      \item[(1)] For every central extension $1\to \TT\to \bar{G}\to
        G\to 1$ the group $\bar{G}$ satisfies BC for $\CC$.
      \item[(2)] $G$ satisfies BC with coefficients in the compact
        operators $\K\cong\K(H)$ for all separable Hilbert spaces $H$
        and with respect to all possible actions of $G$ on $\K$.
      \end{enumerate}
\end{prop}
\begin{proof}
  Assume that (1) holds. Let $\alpha:G\to\Aut(\K)$ be any action of
  $G$ on $\K$ and let $[\om]\in H^2(G,\TT)$ denote the Mackey
  obstruction for this action.  Let
  $$1\to \TT\to G_{\om}\to G\to 1$$
  denote the central extension
  determined by $\om$.  It follows from Lemma \ref{lem-twist} that
  $\alpha$ is Morita equivalent to the twisted action $(\id, \chi_1)$
  of $(G_{\om},\TT)$ on $\CC$.  By assumption, we know that $G_{\om}$
  satisfies BC for $\CC$. It follows from Theorem \ref{thm-ext} that
  $(G_{\om},\TT)$ satisfies BC for $C^*(\TT)\cong C_0(\ZZ)$, or,
  equivalently, that $G$ satisfies BC for $C_0(\ZZ,\K)$ with respect
  to the appropriate action of $G$ (use Proposition \ref{prop-mor}).
  Since our group $G$ does not satisfy directly the assumptions of
  Theorem \ref{thm-ext}, let us briefly explain how it is used: first
  apply part (i) of Theorem \ref{thm-ext} to $N=G_0$, which implies
  that $G$ satisfies BC for $C_0(\ZZ,\K)$ if and only if every compact
  extension $C$ of $G_0$ in $G$ satisfies BC for $C_0(\ZZ,\K)$, and
  then apply part (ii) of Theorem \ref{thm-ext} to the subgroup $\TT$
  of $C_{\om}\subseteq G_{\om}$.
  
  Writing $C_0(\ZZ)=\oplus_{n\in \ZZ}\CC$, the twisted action of
  $(G_{\om},\TT)$ is given by the twisted action $(\id, \chi_n)$ of
  $(G_{\om},\TT)$ on the n'th summand.  Let $q_1:C_0(\ZZ)\to \CC$ be
  the projection on the summand corresponding to $1\in \ZZ$. Consider
  the diagram
  $$
      \begin{CD}
        \k_*^{\top}(G, C_0(\ZZ))@>\mu_{C_0(\ZZ)}>>
        \k_*(C_0(\ZZ)\rtimes_r(G_{\om},\TT))\\
        @Vq_{1,*}VV     @VVq_{1,*}V\\
        \k_*^{\top}(G,\CC) @>>\mu_{\CC}>
        \k_*(\CC\rtimes_r(G_{\om},\TT)).
\end{CD}
$$
(Here the topological $K$-theory $\k_*^{\top}(G,\CC)$ is computed
with respected to the twisted action $(\id,\chi_1)$ of $G\cong
G_{\om}/\TT$ and $\mu_{\CC}$ denotes the twisted assembly map!)  
Since
the vertical arrows are split-surjective and the upper horizontal
arrow is bijective, it follows that the lower horizontal arrow is also
bijective. Thus we see that $(G_{\om},\TT)$ satisfies BC for $\CC$
with respect to the twisted action $(\id, \chi_1)$.  By Morita
equivalence this implies that $G$ satisfies BC for $\K$ with respect
to $\alpha$.

For the opposite direction assume that (2) holds.  Let $1\to \TT\to
\bar{G}\to G\to 1$ be as in (1).  As explained above it follows from
Theorem \ref{thm-ext} that $\bar{G}$ satisfies BC for $\CC$ if
$(\bar{G},\TT)$ satisfies BC for $C^*(\TT)=C_0(\ZZ)$. Using the
stabilization trick, the latter is true if $G$ satisfies BC for
$C_0(\ZZ,\K)$ with respect to an appropriate action of $G$ on
$C_0(\ZZ,\K)$ which fixes the base $\ZZ$.  Using continuity of BC,
this follows easily from the fact that $G$ satisfies BC for arbitrary
actions on $\K$.
   \end{proof}
   
   We also need a result on induced algebras as obtained in
   \cite{CE2}. For this recall that if $H$ is a closed subgroup of $G$
   and $A$ is an $H$-algebra, then the induced algebra $\Ind_H^GA$ is
   defined as
   $$\Ind_H^GA=\{f\in C_b(G,A): f(sh)=h^{-1}(f(s))\;\text{and}\;
   \big(sH\mapsto\|f(s)\|\big)\in C_0(G/H)\}.$$
   Equipped with the
   pointwise operations and the supremum-norm, $\Ind_H^GA$ becomes a
   $C^*$-algebra with $G$-action defined by
   $$s\cdot f(t)=f(s^{-1}t).$$
   The following result follows from
   \cite[Theorem 2.2]{CE2}:

\begin{thm}\label{thm-induced}
  Let $G$, $H$, $A$ and $\Ind_H^GA$ be as above. Then G satisfies BC
  for $\Ind_H^GA$ if and only if $H$ satisfies BC for $A$.
\end{thm}

The result becomes most valuable for us when combined with the
following result of \cite{E}:

\begin{prop}\label{prop-induced}
  Suppose that $H$ is a closed subgroup of $G$ and $B$ is a
  $G$-algebra. Let $\widehat{B}$ denote the set of equivalence classes
  of irreducible representations of $B$ equipped with the usual
  $G$-action defined by $s\cdot\pi(b)=\pi(s^{-1}\cdot b)$.  Then $B$
  is isomorphic (as a $G$-algebra) to $\Ind_H^GA$ for some $H$-algebra
  $A$ if and only if there exists a $G$-equivariant continuous map
  $\varphi:\widehat{B}\to G/H$. Moreover, if $\varphi:\widehat{B}\to
  G/H$ is such a map, then $A$ can be chosen to be $B/I$, with
  $I=\cap\{\ker\pi:\varphi(\pi)=eH\}$ equipped with the obvious
  $H$-action.
\end{prop}

As a corollary of Theorem \ref{thm-induced} and Proposition
\ref{prop-induced} we get in particular:

\begin{cor}\label{cor-induced}
  Suppose that $G$ is a locally compact group and $B$ is a $G$-algebra
  which is type I and such that $G$ acts transitively on
  $\widehat{B}$. Let $\pi\in \widehat{B}$ and let $G_{\pi}$ denote the
  stabilizer of $\pi$ for the action of $G$ on $\widehat{B}$. Then $G$
  satisfies BC for $B$ if and only if $G_{\pi}$ satisfies BC for
  $B/\ker\pi\cong \K(H_{\pi})$, where $H_{\pi}$ denotes the Hilbert
  space of $\pi$.
\end{cor}
\begin{proof}
  Since there is only one orbit for the $G$-action on $\widehat{B}$,
  it follows from results of Glimm \cite{Gli}, that $\widehat{B}$ is
  homeomorphic to $G/G_{\pi}$ via $sG_{\pi}\mapsto s\cdot \pi$. In
  particular, it follows that $\widehat{B}$ is Hausdorff, which
  implies that $B/\ker\pi\cong \pi(B)=\K(H_{\pi})$.  The inverse of
  the above map is clearly a continuous $G$-equivariant map of
  $\widehat{B}$ to $G/G_{\pi}$, and Proposition \ref{prop-induced}
  then implies that $B\cong \Ind_{G_{\pi}}^G(B/\ker\pi)$.  The result
  then follows from Theorem \ref{thm-induced}.
\end{proof}

We now give a short outline of the proof of Theorem \ref{thm-main}.
The main work is required for proving the following proposition:

\begin{prop}\label{prop-main}
  Assume that $G$ is a Lie group with finitely many components and let
  $\alpha:G\to \Aut(\K)$ be an action of $G$ on the compact operators
  on some separable Hilbert space $H$. Then $G$ satisfies BC for $\K$.
\end{prop}

The body of this paper is devoted to give a proof of this proposition
by using induction on the dimension of $G$.  It is fairly easy to see
that the above proposition implies Theorem \ref{thm-main1}. Indeed,
using the first part of Theorem \ref{thm-ext} we can directly reduce
to the case where $G$ is almost connected. Hence Theorem
\ref{thm-main1} follows from

\begin{prop}\label{prop-thm}
  Suppose that Proposition \ref{prop-main} holds. Let $G$ be any
  almost connected group and let $\alpha:G\to \Aut(\K)$ be any action
  of $G$ on the compact operators on some separable Hilbert space $H$.
  Then $G$ satisfies BC with coefficients in $\K$.
\end{prop}
\begin{proof} By the structure theory of almost connected groups
  (e.g. see \cite{MZ}) we can find a compact normal subgroup
  $C\subseteq G$ such that $G/C$ is a Lie group with finitely many
  components.  Using Theorem \ref{thm-ext} we see that $G$ satisfies
  BC for $\K$ if and only if $G/C$ satisfies BC for $\K\rtimes C$
  (with respect to an appropriate twisted action).  Since $C$ is
  compact, it follows that $X:=(\K\rtimes C)\dach$ is discrete, and
  (after stabilizing if necessary), $\K\rtimes C\cong C_0(X,\K)$. Let
  $\tilde{G}:=G/C$ and let $X/\tilde{G}$ denote the space of
  $\tilde{G}$-orbits in $X$.  Since $X$ is discrete, the same is true
  for $X/\tilde{G}$, and we get a decomposition $C_0(X,\K)\cong
  \oplus_{\tilde{G}(x)\in X/\tilde{G}} C_0(\tilde{G}(x),\K)$. By
  continuity of BC (see Proposition \ref{prop-limit}), we conclude
  that $\tilde{G}$ satisfies BC for $C_0(X,\K)$ if and only if
  $\tilde{G}$ satisfies BC for $C_0(\tilde{G}(x),\K)$ for all $x\in
  X$.  Using Corollary \ref{cor-induced}, this will follow if all
  stabilizers $\tilde{G}_x\subseteq \tilde{G}$ satisfy BC for $\K$.
  But since $X$ is discrete, it follows that each stabilizer
  $\tilde{G}_x$ contains the connected component $\tilde{G}_0$ of
  $\tilde{G}$. Thus, each stabilizer is a Lie group with finitely many
  components and the result will follow from Proposition
  \ref{prop-main}.
\end{proof}

As mentioned above, the main idea for the proof of Proposition
\ref{prop-main} is to use induction on the dimension $\dim(G)$ of the
Lie group $G$.  For this we were very much influenced by Puk\'anszky's
proof of the fact that locally algebraic groups (i.e., Lie groups
having the same Lie algebra as some real algebraic group) have type I
group $C^*$-algebras as presented in his recent book \cite{Puk}. We
split the induction argument into two main parts, which deal with the
cases whether $G$ is semi-simple or not. Note that even in the
semi-simple case the result does not follow directly from the existent
results, since all known results only work for trivial coefficients
and require that the groups have finite centers.


\section{Baum-Connes for continuous fields of $C^*$-algebras}
\label{sec-fields}

Let $G$ be a separable locally compact group. Then $G$ is called
$\k$-exact, if the functor $A\mapsto \k_*(A\rtimes_rG)$ is half-exact,
that is: whenever $0\to I\to A\to A/I\to 0$ is a short exact sequence
of $G$-algebras, then the sequence
$$
\k_*(I\rtimes_rG)\to \k_*(A\rtimes_rG)\to \k_*(A/I\rtimes_rG)$$
is
exact in the middle term. Clearly, every exact group is $\k$-exact.
Note that every almost connected group is exact by \cite[Corollary
6.9]{KW1}.

Recall also that an element $\gamma\in \KK_0^G(\CC,\CC)$ is called a
{\em $\gamma$-element} for $G$ if there exists a locally compact
proper $G$-space $Y$, a $C^*$-algebra $D$ equipped with a
nondegenerate and $G$-equivariant $*$-homomorphism $\phi:C_0(Y)\to
ZM(D)$, the center of the multiplier algebra $M(D)$ of $D$, and (Dirac
and dual-Dirac) elements
$$\alpha\in \KK_0^G(D,\CC)\quad \beta\in \KK_0^G(\CC,D)$$
such that
$$\gamma=\beta\otimes_D\alpha\quad\text{and}\quad p_Z^*(\gamma)=1\in
\RKK_0^G(Z; \CC,\CC)$$
for all locally compact proper $G$-spaces $Z$,
where $p_Z:Z\to\{pt\}$.  It is a basic result of Kasparov
\cite[Theorem 5.7]{Kas2} that every almost connected group has a
$\gamma$-element and it follows also from the work of Kasparov (but
see also \cite[\S5]{Tu}) that a $\gamma$-element of $G$ is unique and
that it is an idempotent with the remarkable property that for every
$G$-algebra $B$ the image $\mu_B\big(\k_*^{\top}(G;B)\big)$ of the
assembly map is equal to the $\gamma$-part
$$\gamma\cdot
\k_*(B\rtimes_rG):=\{x\otimes_{B\rtimes_rG}j_G(\sigma_B(\gamma)): x\in
\k_*(B\rtimes_rG)\}.$$
Here and below, we denote by
$j_G:\KK_*^G(A,B)\to \KK_*(A\rtimes_rG, B\rtimes_rG)$ the (reduced)
descent homomorphism of Kasparov and we denote by $\sigma_B:
\KK_*^G(A,D)\to \KK_*^G(B\otimes A, D\otimes B)$ the external tensor
product homomorphism (see \cite[Definition 2.5]{Kas2}).  Note that it
follows from the above discussion that a group $G$ with
$\gamma$-element satisfies BC for a given $G$-algebra $B$ if and only
if $\gamma$ (i.e., $j_G(\sigma_B(\gamma))$) acts as the identity on
$\k_*(B\rtimes_rG)$.  We want to exploit these facts to prove the
following basic result:

\begin{prop}\label{prop-bundle}
  Suppose that $X$ is a separable locally compact space which can be
  realized as the geometric realization of a (probably infinite)
  finite dimensional simplicial complex.  Let $A$ be the algebra of
  $C_0$-sections of a continuous field of $C^*$-algebras $\{A_x: x\in
  X\}$, and let $\alpha:G\to \Aut(A)$ be a $C_0(X)$-linear action of
  $G$ on $A$. Assume further that $G$ is exact and has a
  $\gamma$-element $\gamma\in \KK_0^G(\CC,\CC)$.  Then, if $G$
  satisfies BC with coefficients in each fibre $A_x$, $G$ satisfies BC
  for $A$.
\end{prop}

For the general notion of continuous fields of $C^*$-algebras and
their basic properties we refer to \cite{Fell, ENN, Bl, KW}.

The idea of the proof is to show first that it holds for any closed
interval $I\subseteq \RR$. Then a short induction argument will show
that it holds for any cube in $\RR^n$.  Then the result will follow
from a Mayer-Vietoris argument.  For the proof we first need the
following lemma.

\begin{lem}\label{lem-exact}
  Assume that $G$ is a $\k$-exact group with a $\gamma$-element
  $\gamma\in \KK^G_0(\CC,\CC)$. Let $A$ be a $G$-algebra and let
  $I\subseteq A$ be a $G$-invariant closed ideal of $A$.  Then there
  is a natural six-term exact sequence
  $$
\begin{CD}
  (1-\gamma)\cdot \k_0(I\rtimes_rG) @>>> (1-\gamma)\cdot
  \k_0(A\rtimes_rG)
  @>>> (1-\gamma)\cdot \k_0(A/I\rtimes_rG)\\
  @AAA     @.   @VVV\\
  (1-\gamma)\cdot \k_1(A/I\rtimes_rG) @<<< (1-\gamma)\cdot
  \k_1(A\rtimes_rG) @<<< (1-\gamma)\cdot \k_1(I\rtimes_rG).
\end{CD}
$$
\end{lem}
\begin{proof}
  Since $G$ is $K$-exact, it follows that $A\mapsto \k_*(A\rtimes_rG)$
  is a homotopy invariant and half-exact functor on the category of
  $G$-$C^*$-algebras which also satisfies Bott-periodicity (with
  respect to the trivial $G$-action on $C_0(\RR^2)$).  Then it follows
  from some general arguments (e.g., see \cite[Chapter IX]{B}) that
  there exists a six-term exact sequence
  $$
\begin{CD}
  \k_0(I\rtimes_rG) @>>> \k_0(A\rtimes_rG)
  @>>> \k_0(A/I\rtimes_rG)\\
  @AAA     @.   @VVV\\
  \k_1(A/I\rtimes_rG) @<<< \k_1(A\rtimes_rG) @<<< \k_1(I\rtimes_rG).
\end{CD}
$$
We want to show that all maps in the sequence commute with
multiplication with the $\gamma$-element.  By the construction of the
connecting maps in the above sequence as given in \cite[Chapters VIII
and IX]{B}, it is enough to show that for any pair of $G$-algebras $A$
and $B$ and any $y\in \KK_*^G(A,B)$
$$\k_*(A\rtimes_rG)\to \k_*(B\rtimes_rG); x\mapsto
x\otimes_{A\rtimes_rG}j_G(y)$$
commutes with multiplication with
$\gamma$.  But for this it is enough to show that
$$j_G(y)\otimes_{B\rtimes_rG}j_G(\sigma_B(\gamma))=
j_G(\sigma_A(\gamma))\otimes_{A\rtimes_rG}j_G(y).$$
This follows from
the fact that the descent homomorphism $j_G$ is compatible with
Kasparov products and the fact that
$$y\otimes_B\sigma_B(\gamma)=y\otimes_{\CC}\gamma=\gamma\otimes_{\CC}y
=\sigma_A(\gamma)\otimes_Ay,$$
which follows from \cite[Theorem
2.14]{Kas2}.

It follows now that multiplication with $1-\gamma$ also commutes with
all maps in the above commutative diagram. Since $1-\gamma$ is an
idempotent, it is now easy to see that the full six-term exact
sequence restricts to a six-term exact sequence on the
$1-\gamma$-parts of the respective $\k$-theory groups of the crossed
products.
\end{proof}

\begin{remark}\label{rem-exact}
  It is now a direct consequence of the above proposition that if $G$
  is a $\k$-exact group possessing a $\gamma$-element, and if $0\to
  I\to A\to A/I\to 0$ is a short exact sequence of $G$-algebras, then
  $G$ satisfying BC for two of the algebras in this sequence implies
  that G satisfies BC for all three algebras in the sequence. The same
  result holds without the assumption on the $\gamma$-element (see
  \cite[Proposition 4.2]{CE2} -- which was actually deduced as an easy
  consequence of a result of Kasparov and Skandalis in \cite{KS}).
\end{remark}

We also need the following easy lemma.

\begin{lem}\label{lem-unitary}
  Assume that $X$ is a locally compact space and that $A$ is the
  algebra of $C_0$-sections of the continuous field $\{A_x:x\in X\}$
  of $C^*$-algebras. Assume further that $z\in \k_i(A)$, $i=0,1$, such
  that $q_{x,*}(z)=0$ for some evaluation map $q_x:A\to A_x$. Then
  there exists a compact neighborhood $C$ of $x$ such that
  $q_{C,*}(z)=0$ in $ \k_0(A|_C)$, where $A|_C$ denotes the
  restriction of $A$ to $C$ and $q_C:A\to A|_C$ denotes the quotient
  map.
\end{lem}
\begin{proof} We may assume without loss of generality that $X$ is 
compact.
  Using suspension, it is enough to give a proof for the case $i=0$.
  In what follows, if $B$ is any $C^*$-algebra, we denote by $B^1$ the
  algebra obtained from $B$ by adjoining a unit (even if $B$ is
  already unital). Then $\{A_x^1:x\in X\}$ is a continuous field of
  $C^*$-algebras in a canonical way.  The algebra $\tilde{A}$ of
  sections can be written as the set of pairs $\{(a,f):a\in A, f\in
  C_0(X)\}$ with multiplication given pointwise by the multiplication
  rule of the fibres $A_x^1$. Moreover, we have an obvious unital
  embedding $A^1\to \tilde{A}$.
  
  Assume now that $z\in \k_0(A)$ and $x\in X$ are as in the lemma.  We
  represent $z$ as a formal difference $[p-p']$ for some projections
  $p,p'\in M_l(A^1)$. Since $q_{x,*}(z)=0$ we may assume (after
  increasing dimension if necessary) that there exists a unitary
  $u_x\in M_l(A_x^1)$ such that $u_xp_xu_x^*=p'_x$. After passing to
  $\left(\begin{smallmatrix} u&0\\ 0&u^*\end{smallmatrix}\right)$ if
  necessary, we may further assume that $u_x$ lies in the connected
  component of the identity of $U(M_l(A_x^1))$.  Thus, there exists a
  unitary $u\in M_l(A^1)$ such that $q_x(u)=u_x$. Since $u$ is a
  continuous section in $\tilde{A}$, it follows that there exists a
  compact neighborhood $C$ of $x$ such that $\|u_yp_yu_y^*-p'_y\|< 1$
  for all $y\in C$, which implies that $[p_C]=[u_Cp_Cu_C^*]=[p'_C]\in
  \k_0(A|_C^1)$, where $p_C, u_C,$ and $p'_C$ denote the restrictions
  of $p,u,p'$ to $C$, respectively. But this shows that
  $q_{C,*}(z)=[p_C-p'_C]=0$ in $ \k_0(A|_C)$.
\end{proof}

\begin{proof}[Proof of Proposition \ref{prop-bundle}]
  Since $G$ is exact, it follows from \cite[Theorem]{KW} that the
  crossed products $\{A_x\rtimes_rG: x\in X\}$ form a continuous
  bundle such that $A\rtimes_rG$ is the algebra of continuous sections
  of this bundle.  We start with proving the result in the special
  case where $X=[0,1]\subseteq \RR$.
  
  Recall from the above discussions that $G$ satisfies BC for a given
  $G$-algebra $B$ if and only if $(1-\gamma)\cdot
  \k_*(B\rtimes_rG)=\{0\}$.  In particular, it follows from our
  assumptions that $(1-\gamma)\cdot \k_*(A_x\rtimes_rG)=\{0\}$ for all
  $x\in I$.  Assume now that $z\in (1-\gamma)\cdot \k_i(A\rtimes_rG)$,
  $i=0,1$, and let $q_x:A\rtimes_rG\to A_x\rtimes_rG$ denote the
  evaluation maps for each $x\in X$.  Then $q_{x,*}(z)\in
  (1-\gamma)\cdot \k_i(A_x\rtimes_rG)=\{0\}$ for all $x\in I$.  Thus,
  using Lemma \ref{lem-unitary}, we see that there exists a partition
  $0=x_0<x_1<\cdots <x_l=1$ such that $q_{[x_{j-1},x_j],*}(z)=0$ in $
  \k_i(A|_{[x_{j-1},x_j]}\rtimes_rG)$.  Now let
  $O=[0,1]\setminus\{x_0,\dots, x_l\}$ and let $A|_O=C_0(O)\cdot
  A\cong \bigoplus_{j=1}^l A|_{(x_{j-1},x_j)}$.  It follows from the
  exact sequence
  $$
  (1-\gamma)\cdot \k_i(A_O\rtimes_rG) \to (1-\gamma)\cdot
  \k_i(A\rtimes_rG) \to \bigoplus_{j=0}^l (1-\gamma)\cdot
  \k_i(A_{x_j}\rtimes_rG)=\{0\}
  $$
  that there exists a $z'\in (1-\gamma)\cdot \k_i(A_O\rtimes_rG)$
  such that $z$ is the image of $z'$ under the inclusion.  Since
  $$(1-\gamma)\cdot
  \k_i(A_O\rtimes_rG)=\bigoplus_{j=1}^l(1-\gamma)\cdot
  \k_i\big(A|_{(x_{j-1},x_j)}\rtimes_rG\big),$$
  we may write
  $z'=\sum_{j=1}^l z'_j$ with $z'_j\in (1-\gamma)\cdot
  \k_i\big(A|_{(x_{j_1},x_j)}\rtimes_rG\big)$ for each $1\leq j\leq
  l$.  Thus it is enough to show that $z'_j=0$ for each $1\leq j\leq
  l$.  In what follows, we write $A_j=A|_{(x_{j-1}, x_j)}$ and
  $\bar{A}_j=A|_{[x_{j-1}, x_j]}$.  Since $(1-\gamma)\cdot
  \k_i(A_{x_k}\rtimes_rG)=\{0\}$ for all $0\leq k\leq l$ we obtain a
  six-term exact sequence
  $$
\begin{CD}
  (1-\gamma)\cdot \k_0(A_j\rtimes_rG) @>>> (1-\gamma)\cdot
  \k_0(\bar{A}_j\rtimes_rG)
  @>>> 0\\
  @AAA     @.   @VVV\\
  0 @<<< (1-\gamma)\cdot \k_1(\bar{A}_j\rtimes_rG) @<<<
  (1-\gamma)\cdot \k_1(A_j\rtimes_rG).
\end{CD}
$$
Since the image of $z'_j$ in $\k_i(\bar{A}_j\rtimes_rG)$ coincides
with the image of $z$ in $\k_i(\bar{A}_j\rtimes_rG)$, we see that
$z_j'$ maps to $0$ under the isomorphism $(1-\gamma)\cdot
\k_0(A_j\rtimes_rG)\to (1-\gamma)\cdot \k_0(\bar{A}_j\rtimes_rG)$, so
$z'_j=0$.

We now show by induction on $n$ that the result is true for
$[0,1]^n\subseteq \RR^n$.  For this assume that $\{A_x: x\in
[0,1]^n\}$ is a continuous field over the cube and $A$ is the algebra
of continuous sections of this field. We write $[0,1]^n=\cup_{y\in
  [0,1]}\{y\}\times [0,1]^{n-1}$ and put $A_y=A|_{\{y\}\times
  [0,1]^{n-1}}$.  Then $\{A_y:y\in [0,1]\}$ is a continuous field over
$[0,1]$ and $A$ is also the section algebra of this bundle. If
$\alpha$ is a $C([0,1]^n)$-linear action on $A$, it is also
$C([0,1])$-linear with respect to the bundle structure of $A$ over
$[0,1]$ coming from the above decomposition of the cube.  Moreover,
the actions on the fibres $A_y$ are clearly $C([0,1]^{n-1})$-linear,
so by the induction assumption we know that $G$ satisfies BC with
coefficients in $A_y$ for all $y\in [0,1]$. We now apply the above
result to the bundle $\{A_y:y\in [0,1]\}$ to conclude that $G$
satisfies BC with coefficients in $A$.

In a next step we show that the result holds for the open cubes
$(0,1)^n\subseteq \RR^n$. By similar arguments as given above it
suffices to show that the result holds for open intervals.  So assume
that $\{A_x:x\in (0,1)\}$ is a continuous field with section algebra
$A$. Let $x_1<x_2\in (0,1)$.  Then it follows from the first part of
the proof that $G$ satisfies BC with coefficients in $A_{[x_1,x_2]}$.
Since, by assumption, $G$ also satisfies BC for the fibres, a six-term
sequence argument shows that it also satisfies BC with coefficients in
$(x_1, x_2)$.  Writing $A=\lim_{n\to\infty}A|_{(\frac{1}{n},
  1-\frac{1}{n})}$ and using continuity of the BC conjecture, it
follows that $G$ satisfies BC for $A$.

Since the result of the proposition is clearly invariant under
replacing the space $X$ by a homeomorphic space $Y$, we now see that
the result holds for all open or closed simplices. We now proof the
general result for simplicial complexes via induction on the dimension
of the complex.  By continuity of the conjecture, the result is clear
for zero-dimensional complexes. If $X$ has dimension $n$, let $W_n$
denote the interiors of all $n$-dimensional simplices in $X$. Then
$W_n$ is homeomorphic to a disjoint union of open $n$-dimensional
cubes, so the result holds for $W_n$. Since $X\smallsetminus W_n$ is a
simplicial complex of dimension $n-1$, the result is true for
$X\smallsetminus W_n$ by the induction assumption. The result then
follows from another easy application of the six-term sequence (see
Remark \ref{rem-exact}).
\end{proof}

\begin{remark}\label{rem-countable}
  Let $G$ be a second countable locally compact group and let $X$ be a
  second countable locally compact $G$-space. Following Glimm we say
  that the quotient space $X/G$ is {\em countably separated}, if all
  orbits $G(x)$ are locally closed, i.e. $G(x)$ is open in its closure
  $\overline{G(x)}$.  Glimm showed in \cite[Theorem]{Gli} that $X/G$
  being countably separated is equivalent to each of the following
  conditions:
\begin{itemize}
\item The canonical map $G/G_x\to G(x), gG_x\mapsto g\cdot x$ is a
  homeomorphism for each $x\in X$.
\item There exists a sequence of $G$-invariant open subsets
  $\{U_{\nu}\}_{\nu}$ of $X$, where $\nu$ runs through the ordinal
  numbers such that
      \begin{enumerate}
      \item $U_{\nu}\subseteq U_{\nu+1}$ for each $\nu$ and
        $\big(U_{\nu+1}\smallsetminus U_{\nu})/G$ is Hausdorff.
      \item If $\nu$ is a limit ordinal, then
        $U_{\nu}=\cup_{\mu<\nu}U_{\mu}$.
      \item There exits an ordinal number $\nu_0$ such that
        $X=U_{\nu_0}$.
\end{enumerate}
\end{itemize}

Unfortunately, in order to apply our bundle results to such actions,
we need to know that the spaces $\big(U_{\nu+1}\smallsetminus
U_{\nu}\big)/G$ have a simplicial structure as required in Proposition
\ref{prop-bundle}.  However, we can summarize the above results to
prove the following theorem which turns out to be sufficient for our
purposes:
\end{remark}

\begin{thm}\label{thm-BCbundle}
  Suppose that $G$ is a second countable exact group which has a
  $\gamma$-element. Let $X$ be a second countable locally compact
  $G$-space, let $\{A_x:x\in X\}$ be a continuous bundle of
  $C^*$-algebras over $X$ with section algebra $A$ and let
  $\alpha:G\to \Aut(A)$ be an action of $G$ on $A$ which is compatible
  with the given action of $G$ on $X$. Assume further that the
  following assumptions are satisfied:
      \begin{enumerate}
      \item $X/G$ is countably separated.
      \item There exists an increasing sequence of open $G$-invariant
        subsets $\{U_{\nu}\}_{\nu}$ of $X$ such that
        $U_{\nu}=\cup_{\mu<\nu}U_{\mu}$ if $\nu$ is a limit ordinal
        and $X=U_{\nu_0}$ for some $\nu_0$.
      \item For each $\nu$ there exists a finite dimensional
        simplicial complex $Y_{\nu}$ and a continuous and open
        surjection $q_{\nu}:U_{\nu+1}\smallsetminus U_{\nu}\to
        Y_{\nu}$ such that for all $y\in Y_{\nu}$,
        $q_{\nu}^{-1}(\{y\})$ is a finite union of orbits in $X$.
      \item Each stabilizer $G_x$, $x\in X$, satisfies BC for $A_x$.
     \end{enumerate}
     Then $G$ satisfies BC for $A$.
\end{thm}
\begin{proof} For each ordinal $\nu$ let $A_{\nu}:=C_0(U_{\nu})A$
  denote the ideal corresponding to $U_{\nu}$. We show by transfinite
  induction that $G$ satisfies BC with coefficients in $A_{\nu}$ for
  each $\nu$.  Since $A=A_{\nu_0}$ for some $\nu_0$, the result will
  follow.
  
  We start by showing that $G$ satisfies BC with coefficients in
  $A_1$.  Since the map $q_1:U_1\to Y_1$ is open, we can regard $A_1$
  as a section algebra of a continuous bundle over $Y_1$ with fibres
  isomorphic to $A|_{q_1^{-1}(y)}$.  By Proposition \ref{prop-bundle}
  it is therefore sufficient to prove that $G$ satisfies BC for
  $A|_{q_1^{-1}(y)}$ for all $y\in Y_1$.  Fix $y\in Y_1$ and put
  $Z:=q_1^{-1}(y)$. Since $Z$ is a finite union of $G$-orbits, we find
  a finite sequence
  $$Z=Z_0\supseteq Z_1\supseteq \cdots \supseteq Z_l=\emptyset$$
  of
  open invariant subsets of $Z$ such that $\big(Z_{i-1}\smallsetminus
  Z_{i}\big)/G$ is a discrete finite set. To see this let $C_1$ be the
  union of all closed $G$-orbits in $Z$ (such orbits must exist by the
  finiteness of $Z/G$!).  Then $C_1$ is closed in $Z$ and $C_1/G$ is
  discrete. Put $Z_1=Z\smallsetminus C_1$ and then define the $Z_i$'s,
  $i>1$, inductively by the same procedure.  Using six-term sequences
  (e.g., see Remark \ref{rem-exact}), $G$ satisfies BC for $A|_Z$ if
  $G$ satisfies BC for all $A|_{Z_{i-1}\smallsetminus Z_i}$, which in
  turn follows if $G$ satisfies BC for $A|_{G(x)}$ for any $G$-orbit
  $G(x)\subseteq X$.  Since $G(x)$ is homeomorphic to $G/G_x$ via the
  canonical map (see Remark \ref{rem-countable} above), it follows
  from Proposition \ref{prop-induced}, that $A|_{G(x)}$ is
  $G$-equivariantly isomorphic to $\Ind_{G_x}^GA_x$. Since, by
  assumption, $G_x$ satisfies BC with coefficients in $A_x$, it
  follows then from Theorem \ref{thm-induced} that $G$ satisfies BC
  with coefficients in $A|_{G(x)}$.  This completes the proof for
  $A_1$.
  
  Assume now that $\nu$ is an ordinal number and that we have already
  shown that $G$ satisfies BC for $A_{\mu}$ for all $\mu<\nu$.  If
  $\nu=\mu+1$ for some ordinal $\mu$, it follows from the same
  reasoning as for the case $\nu=1$ that $G$ satisfies BC for
  $A|_{U_{\nu}\smallsetminus U_{\mu}}\cong A_{\nu}/A_{\mu}$. Since $G$
  satisfies BC for $A_{\mu}$ by the induction assumption, it follows
  from Remark \ref{rem-exact} that $G$ satisfies BC for $A_{\nu}$.
  
  Assume now that $\nu$ is a limit ordinal and $G$ satisfies BC for
  each $\mu<\nu$. Then $U_{\nu}=\cup_{\mu<\nu}U_{\mu}$ which implies
  that $A_{\nu}=\lim_{\mu<\nu}A_{\mu}$ is the inductive limit of the
  $A_{\mu}$.  Thus it follows from Proposition \ref{prop-limit} that
  $G$ satisfies BC for $A_{\nu}$.
\end{proof}

\section{The semi-simple case}\label{sec-semi}

In this section we want to show that Proposition \ref{prop-main} is
true if $G$ is semi-simple.  For this we first have to obtain a slight
extension of Lafforgue's results on the Baum-Connes conjecture for
semi-simple groups with finite center.

Let us first recall the basic idea of Lafforgue's proof of the
Baum-Connes conjecture for such groups. If $G$ is a locally compact
group we let $C_c(G)$ denote the convolution algebra of $G$ consisting
of continuous functions with compact supports.  A norm $\|\cdot\|$ on
$C_c(G)$ is called {\em good} if convolution is continuous with
respect to this norm and if $\|f\|$ only depends on the absolute value
of $f$ for all $f\in C_c(G)$ (i.e., $\|f\|=\| |f| \|$ for all $f\in
C_c(G)$).  A {\em good completion} $\mathcal A(G)$ of $C_c(G)$ is a
completion with respect to a good norm on $C_c(G)$.  Note that
$L^1(G)$ is always a good completion of $C_c(G)$, but $C^*(G)$ and
$C_r^*(G)$ are in general not good completions of $C_c(G)$.

If $\mathcal A(G)$ is a good completion of $C_c(G)$, then Lafforgue
constructed an assembly map
$$\mu_{\mathcal{A}(G)}: \k_*^{\top}(G,\CC)\to \k_*(\mathcal A(G)).$$
Moreover, if the identity on $C_c(G)$ extends to a continuous
embedding $\iota: \mathcal A(G)\to C_r^*(G)$, he also shows that the
assembly map $\mu: \k_*^{\top}(G,\CC)\to \k_*(C_r^*(G))$ factors
through $ \k_*(\mathcal A(G))$, i.e.,
$$\mu=\iota_*\circ \mu_{\mathcal A(G)}$$
(see \cite[Proposition
1.7.6]{Laf}).  Thus, if we know that $\mu_{\mathcal{A}(G)}$ is an
isomorphism for all good completions of $C_c(G)$, and if we further
know that there exists a good completion $\mathcal A(G)\subseteq
C_r^*(G)$ such that the inclusion $\iota_*: \k_*(\mathcal A(G))\to
\k_*(C_r^*(G))$ is an isomorphism, it follows that $\mu:
\k_*^{\top}(G,\CC)\to \k_*(C_r^*(G))$ is an isomorphism.  Note that
$\iota_*:\k_*(\mathcal A(G))\to \k_*(C_r^*(G))$ is an isomorphism
whenever $\mathcal A(G)$ is closed under holomorphic functional
calculus in $C_r^*(G)$. Now Lafforgue was able to prove the following
deep results:

\begin{thm}[{cf \cite[Th\'eor\`eme 17.13 and Chapitre 3]{Laf}}]
\label{thm-Laf1}
  Assume that $G$ is a second countable locally compact group such
  that $G$ acts isometrically and properly on a Riemannian manifold
  with nonpositive sectional curvature which is bounded below.  Then
  $\mu_{\mathcal A(G)}: \k_*^{\top}(G,\CC)\to \k_*(\mathcal A(G))$ is
  an isomorphism for every good completion $\mathcal A(G)$ of
  $C_c(G)$.
\end{thm}

In fact, Lafforgue was even able to show that the above result holds
with arbitrary $C^*$-algebra coefficients, but we do not need this
more general result here.  Note that if $G$ is semi-simple with finite
center, then the Riemannian manifold of the theorem can be chosen to
be the symmetric space $G/K$, where $K$ is the maximal compact
subgroup of $G$.

In the second step for the proof of BC for semi-simple groups,
Lafforgue constructed a {\em Schwartz-algebra} $\mathcal S(G)\subseteq
C_r^*(G)$ which is a good completion of $C_c(G)$ which is closed under
holomorphic functional calculus in $C_r^*(G)$. In fact, this
construction followed a more general principle, which we are now going
to describe in more detail.

Assume that $G$ is a unimodular group and $K$ is a compact subgroup of
$G$. Then, following Lafforgue (see \cite[Chapitre 4]{Laf}), we say
that the pair $(G,K)$ satisfies property $(HC)$ if the following
conditions are satisfied
\begin{enumerate}
\item[(HC1)] There exists a continuous function $d:G\to [0,\infty)$
  such that $d(e)=0$, $d(kgk')=d(g)$ for all $k,k'\in K, g\in G$ and
  $d(gg')\leq d(g)+d(g')$ for all $g,g'\in G$.
\item[(HC2)] There exists a continuous function $\phi:G\to (0,1]$ such
  that $\phi(e)=1$, $\phi(g^{-1})=\phi(g)$, $\phi(kgk')=\phi(g)$ for
  all $g\in G$ and $k,k'\in K$, and
  $$\int_K\phi(gkg')\,dk=\phi(g)\phi(g')\quad\text{for all $g,g'\in
    G$},$$
  with respect to the normalized Haar measure on $K$.
\item[(HC3)] There exists a $t_0\in \RR$ such that
  $\big(t\mapsto\phi(g)(1+d(g))^{-t}\big)\in L^2(G)$ for all $t>t_0$.
\end{enumerate}
If $G$ is a connected semi-simple group with finite center, and if $K$
is the maximal compact subgroup of $G$, then Lafforgue was able to
show that $(G,K)$ satisfies (HC). The following theorem then completes
the proof of BC for connected semi-simple groups with finite center.

\begin{thm}[{cf \cite[Proposition 4.2.1]{Laf}}]\label{thm-Laf2}
  Assume that $G$ is a unimodular group and $K$ is a compact subgroup
  of $G$ such that the pair $(G,K)$ satisfies (HC).  Then there exists
  a good completion $\mathcal S(G)\subseteq C_r^*(G)$ of $C_c(G)$ such
  that $\mathcal S(G)$ is closed under holomorphic functional calculus
  in $C_r^*(G)$.
\end{thm}

In order to extend Lafforgue's methods to extensions of semi-simple
groups, we first show that property (HC) is closed under compact
extensions.

\begin{lem}\label{lem-below}
  Assume that $(G,K)$ satisfies property (HC). Assume further that
  $$1\to C\to \tilde{G}\stackrel{q}{\to} G\to 1$$
  is a group extension
  with $C$ compact. Let $\tilde{K}=q^{-1}(K)\subseteq \tilde{G}$.
  Then $(\tilde{G},\tilde{K})$ satisfies (HC).
\end{lem}
\begin{proof} Since compact extensions of unimodular groups
  are unimodular, $\tilde{G}$ is unimodular.  Let $(d,\phi)$ be a pair
  of functions which satisfy HC1, HC2, and HC3 with respect to
  $(G,K)$. Define $\tilde{d}(g)=d(q(g))$ and
  $\tilde\phi(g)=\phi(q(g))$. Then a straightforward computation shows
  that $(\tilde{d},\tilde{\phi})$ satisfies HC1, HC2, and HC3 with
  respect to $(\tilde{G},\tilde{K})$.
\end{proof}

The second result is slightly more technical.

\begin{lem}\label{lem-upper}
  Assume that $G$ is a unimodular Lie group with finitely many
  components and let $K$ be a maximal compact subgroup of $G$.  Let
  $G_0$ denote the connected component of $G$ and let $K_0$ be a
  maximal compact subgroup of $G_0$ such that $(G_0, K_0)$ satisfies
  (HC).  Then $(G,K)$ satisfies (HC), too.
\end{lem}
\begin{proof}
  First note that we may assume that $K_0=K\cap G_0$.  To see this
  observe first that, since $K_0$ is a compact subgroup of $G$, we may
  assume without loss of generality that $K_0\subseteq K\cap G_0$.
  But the maximality of $K_0$ then implies equality.  It follows from
  this that $K_0$ is a normal subgroup of $K$ and since $G/K$ is
  connected, it follows that the inclusion $K\to G$ induces a group
  isomorphism $K/K_0\cong G/G_0$.

  Let $(d_0,\phi_0)$ be a pair of functions satisfying conditions
  HC1--HC3 for $(G_0,K_0)$. It follows from the above remarks that we
  can write every element of $G$ as a product $kg$ with $k\in K$,
  $g\in G_0$.  We then define
  $$d(kg)=\int_K d_0(lgl^{-1})\, dl\quad\text{and}\quad
  \phi(kg)=\int_K \phi_0(lgl^{-1})\,dl.$$
  To see that $d$ and $\phi$
  are well defined assume that we have two factorizations $kg=k'g'$
  with $k,k'\in K, g,g'\in G_0$.  Then $g=k^{-1}k'g'$ with
  $h:=k^{-1}k'\in K_0$.  Since $K$ normalizes $K_0$, it follows from
  the left and right $K_0$-invariance of $d_0$ that
      \begin{align*}
        d(kg)&=\int_K d_0(lgl^{-1})\,dl=\int_K d_0(lhg'l^{-1})\,dl
        =\int_K d_0\big((lhl^{-1})(lg'l^{-1})\big)\,dl\\
        &=\int_K d_0(lg'l^{-1})\,dl =d(k'g').
      \end{align*}
      A similar computation shows that $\phi$ is well defined.
      
      We are now going to check properties HC1--HC3 for $(d,\phi)$.
      It follows directly from the definition of $d$ and $\phi$ that
      they are left invariant under the action of $K$. To see right
      invariance, we compute for $h\in K$:
      $$
      d(kgh)=d(khh^{-1}gh)=d(h^{-1}gh) =\int_K
      d_0(lh^{-1}ghl^{-1})\,dl \stackrel{l\mapsto lh}{=}\int_K
      d_0(lhl^{-1})\,dl=d(kg).
      $$
      So $d$ is also right invariant and a similar computation show
      that the same is true for $\phi$.  Since Haar measure on $K$ is
      normalized, it follows that $d(e)=0$ and $\phi(e)=1$. Moreover,
      if $kg, hg'\in G$ with $k,h\in K, g,g'\in G_0$ we get
      \begin{align*}
        d(kghg')&=d(khh^{-1}ghg')=\int_K d_0(lh^{-1}ghg'l^{-1})\,dl
        =\int_K d_0(lh^{-1}ghl^{-1}lg'l^{-1})\,dl\\
        &\leq \int_K d_0(lh^{-1}ghl^{-1})\,dl +\int_K
        d_0(lg'l^{-1})\,dl =d(kg)+d(hg').
       \end{align*}
       In order to prove the multiplication rule for $\phi$ we use
       Weil's formula
       $$\int_K \varphi(l)\,dl=
       \int_{K/K_0}\left(\int_{K_0}\varphi(lm)\,dm\right)\,d\dot{l}
       =\int_K\left(\int_{K_0}\varphi(lm)\,dm\right)\,dl,$$
       with
       respect to normalized Haar measures on $K, K_0$, and $K/K_0$ to
       compute
       \begin{align*}
         \int_K \phi(kg l hg')\,dl&{=}\int_K\phi(gl hg'h^{-1})\,dl
         =\int_K\int_{K_0}\phi(glmhg'h^{-1})\,dm\,dl\\
         &\stackrel{(*)}{=}\int_K\int_{K_0}\phi(l^{-1}glmhg'h^{-1})\,dm\,dl\\
         &=\int_K\int_{K_0}\int_K\phi_0(nl^{-1}glmhg'h^{-1}n^{-1})
\,dn\,dm\,dl         \\
         &\stackrel{(**)}{=}
         \int_K\int_K\int_{K_0}\phi_0(nl^{-1}gln^{-1}mnhg'h^{-1}n^{-1})
\,dm\,dl\,dn \\
         &\stackrel{(***)}{=}
         \int_K\int_K\phi_0(nl^{-1}gln^{-1})\phi_0(nhg'h^{-1}n^{-1})
\,dl\,dn\\
         &\stackrel{l\mapsto l^{-1}n}{=}
         \left(\int_K\phi_0(lgl^{-1})\,dl\right)
         \left(\int_K\phi_0(ng'n^{-1})\,dn\right)\\
         &=\phi(kg)\phi(hg').
      \end{align*}
      Here the equation (*) follows from the $K$-invariance of $\phi$,
      (**) follows from Fubini together with the transformation
      $m\mapsto n^{-1}mn$, and (***) follows from property HC2 for
      $\phi_0$.  This completes the proof of HC1 and HC2.
      
      For the proof of HC3 we write $F_s(g):= \phi(g)(1+d(g))^{-s}$,
      $s>0$. Since $G=KG_0$, it follows from the $K$-invariance of
      $\phi$ and $d$ that the integrals of $F_s^2$ over the
      $G_0$-cosets coincide. Thus, since $G/G_0$ is finite, it is
      enough to show that $F_s|_{G_0}\in L^2(G_0)$ for some $s>0$. For
      this we first choose a set of representatives $t_1,\ldots,
      t_n\in K$ for $K/K_0$ with $t_1=e$.  Then, for $g\in G_0$, we
      obtain the inequality
      $$
      1+d(g)=1+\frac{1}{n}\sum_{i=1}^n d_0(t_igt_i^{-1})\geq
      1+\frac{1}{n}d_0(g),
      $$
      from which it follows that $(1+d(g))^n\geq 1+d_0(g)$ for all
      $g\in G_0$.  Thus, if $t\in \RR$ such that $\left(g\mapsto
        \phi_0(g)(1+d_0(g))^{-t}\right)\in L^2(G_0)$, then we also
      have $\left(g\mapsto \phi_0(g)(1+d(g))^{-nt}\right)\in
      L^2(G_0)$.  So let $s=nt$ with $t$ as above. Then we get
\begin{align*}
  F_s(g)&=\phi(g)(1+d(g))^{-s}=
  \frac{1}{n}\sum_{i=1}^n\phi_0(t_igt_i^{-1})(1+d(g))^{-s}\\
  &=\frac{1}{n}\sum_{i=1}^n\phi_0(t_igt_i^{-1})(1+d(t_igt_i^{-1}))^{-s}.
\end{align*}
Since $G$ is unimodular, it follows that each summand $\left(g\mapsto
  \phi_0(t_igt_i^{-1})(1+d(t_igt_i^{-1}))^{-s}\right)\in L^2(G_0)$,
and hence $F_s|_{G_0}\in L^2(G_0)$.
\end{proof}

We are now ready to combine the above results to get

\begin{prop}\label{prop-BC}
  Let $G$ be a locally compact group with finitely many components.
  Assume further that $G$ has a compact normal subgroup $C\subseteq
  G_0$ such that $G_0/C$ is semi-simple with finite center. Then $G$
  satisfies BC with trivial coefficients.
\end{prop}
\begin{proof} Let $K_0$ denote the maximal compact subgroup of $G_0$.
  Then $G_0/K_0\cong (G_0/C)/(K_0/C)$ is a symmetric space and
  therefore has nonpositive sectional curvature. Moreover, if $K$ is a
  maximal compact subgroup of $G$ such that $K\cap G_0=K_0$, we see
  that $G/K\cong G_0/K_0$ as a Riemannian manifold.  Since $G$ acts
  isometrically and properly on $G/K$, $G$ satisfies the assumptions
  of Theorem \ref{thm-Laf1}. By Lafforgue's results we also know that
  $(G_0/C,K_0/C)$ satisfies (HC).  Lemmas \ref{lem-below} and
  \ref{lem-upper} then imply that $(G,K)$ also satisfies property
  (HC). Thus, it follows from the combination of Theorem
  \ref{thm-Laf1} with Theorem \ref{thm-Laf2} that $G$ satisfies BC
  with coefficients in $\CC$.
\end{proof}

Using the results on continuous fields of actions as presented in the
previous section, we are now able to prove

\begin{prop}\label{prop-reduct}
  Assume that $G$ is a Lie group with finitely many components such
  that $G_0$ is reductive, i.e., the Lie algebra $\frak{g}$ of $G$ is
  a direct sum of two ideals $\frak{g}=\frak{s}\oplus\frak{z}$ with
  $\frak{s}$ semi-simple and $\frak{z}$ abelian.  Then $G$ satisfies
  BC for $\CC$.
\end{prop}
\begin{proof}
  Let $Z=Z(G_0)$ denote the center of $G_0$. Using Theorem
  \ref{thm-ext} it is enough to show that $G/Z$ satisfies BC with
  coefficients in $C^*(Z)\otimes\K\cong C_0(\widehat{Z},\K)$, where
  the action of $G/Z$ on the dual space $\widehat{Z}$ of $Z$ is given
  via conjugation. Since $Z$ is central in $G_0$, it follows that this
  action factors through an action of the finite group $G/G_0$.
  Moreover, since $\widehat{Z}$ is a manifold (since $Z$ is a
  compactly generated abelian group), it follows that the quotients of
  the orbit-types in $\widehat{Z}$ are manifolds. From this we easily
  obtain a finite decomposition sequence
  $$\emptyset=U_0\subseteq U_1\subseteq\cdots\subseteq
  U_l=\widehat{Z}$$
  of open $G$-invariant subsets of $\widehat{Z}$
  such that the quotients of the differences $U_j\smallsetminus
  U_{j-1}$ are homeomorphic to geometric realizations of finite
  dimensional simplicial complexes.  Moreover, since all stabilizers
  for the action of $G/Z$ on $\widehat{Z}$ contain $G_0/Z$, which is
  semi-simple with trivial center, it follows from a combination of
  Proposition \ref{prop-BC} with Proposition \ref{prop-central} that
  all stabilizers satisfy BC for $\K$. The result then follows from
  Theorem \ref{thm-BCbundle}.
\end{proof}

Since any central extension of a semi-simple group is reductive, we
now get the desired result for general semi-simple groups.

\begin{cor}\label{cor-semi-simple}
  Let $G$ be a semi-simple Lie group with finitely many components and
  let
  $$1\to \TT\to \bar{G}\to G\to 1$$
  be a central extension of $G$ by
  $\TT$.  Then $\bar{G}$ satisfies BC for $\CC$.  As a consequence
  (using Proposition \ref{prop-central}), $G$ satisfies BC for $\K$
  with respect to arbitrary actions of $G$ on $\K$.
\end{cor}

\section{The general case}\label{sec-perfect}

We now want to give a proof of Proposition \ref{prop-main}. As
indicated in the first section, we are going to use an induction
argument on the dimension $n=\dim(G)$.  Since any one-dimensional Lie
group with finitely many components is amenable, and since amenable
groups satisfy BC for arbitrary coefficients, the case $n=1$ is clear.
Assume now that $G$ is an arbitrary Lie group with finitely many
components. Let $G_0$ denote the connected component of $G$.  Let $N$
denote the nilradical of $G$ and let $\frak{n}$ and $\frak{g}$ denote
the Lie algebras of $N$ and $G$, respectively.  If $\frak{n}=\{0\}$,
then $G$ is semi-simple and the result follows from the previous
section. So we may assume that $\frak{n}\neq\{0\}$.

It is shown in \cite[Lemma 4 on p. 24]{Puk} that the subgroup $H$ of
$G_0$ corresponding to the subalgebra
$\frak{h}=[\frak{g},\frak{g}]+\frak{n}$ of $\frak{g}$ is closed in
$G_0$.  Further, if $\frak{s}$ is a Levi section in $\frak{g}$, i.e.,
$\frak{s}$ is a maximal semi-simple subalgebra of $\frak{g}$, then
$\frak{h}=\frak{s}+\frak{n}$ (e.g., see the discussion in the proof of
\cite[Sublemma on p. 24]{Puk}). In particular, $H/N$ is semi-simple.
Clearly, $H$ is a normal subgroup of $G$ and $G/H$ is a finite
extension of a connected abelian Lie group.  Let $M$ denote the
inverse image of the maximal compact subgroup of $G/H$ in $G$. It
follows then from Theorem \ref{thm-ext} that $G$ satisfies BC for $\K$
if and only if $M$ satisfies BC for $\K$. Note that the connected
component of $M/H$ is a compact connected abelian Lie group, hence a
torus group.

Thus, replacing $G$ by $M$, we may from now on assume that $G$ has the
following structure: There exist closed normal subgroups
\begin{equation}\label{eq-structure}
      N\subseteq H\subseteq G_0\subseteq G
\end{equation}
such that $N$ is a non-trivial connected nilpotent Lie group, $H/N$ is
semi-simple, $G_0/H$ is a torus group and $G/G_0$ is finite. Moreover,
by induction we may assume that every almost connected Lie group with
smaller dimension satisfies BC for $\K$ with respect to arbitrary
actions on $\K$, or, equivalently (by Proposition \ref{prop-central}),
every central extension by $\TT$ satisfies BC for $\CC$.
It is now useful to recall the following result of Chevalley
(see \cite[Proposition 5, p.324]{Ch}):

\begin{prop}\label{prop-liealg} Let $\frak{g}\subseteq \frak{gl}(V)$
be a Lie-algebra of endomorphisms of the finite dimensional
real vector space $V$. Then 
$\frak{g}$ is algebraic (i.e. it corresponds to a real 
algebraic subgroup $G\subseteq \GL(V)$)
if and only if there exist subalgebras $\frak{s},\frak{a}, \frak{n}$
of $\frak{g}$ with $\frak{g}=\frak{s}+\frak{a}+ \frak{n}$,
$\frak{s}$ is semi-simple, $\frak{n}$ is the largest ideal 
of $\frak{g}$ consisting of nilpotent
endomorphisms, and $\frak{a}$ is an algebraic abelian 
subalgebra of $\frak{gl}(V)$ consisting of semi-simple
endomorphisms such that $[\frak{s},\frak{a}]
\subseteq \frak{a}$.
\end{prop}

Using Ado's theorem
(see \cite[Th\'eor\`eme 5 on p. 333]{Ch})
and Proposition \ref{prop-liealg}, it follows that the group $H$
considered above is locally algebraic,
i.e., the Lie algebra $\frak{h}$ has a faithful representation 
as an algebraic Lie subalgebra into some $\frak{gl}(V)$.

Using this structure, the main idea is to apply the Mackey machine to
a suitable abelian subgroup $S$ of $N$ which is normal in $G$.  The
fact that $G$ is very close to an algebraic group implies that the
action of $G$ on the dual $\widehat{S}$ of $S$ has very good
topological properties, which is precisely what we need to make
everything work. As a first hint that this approach is feasible we
prove:

\begin{lem}\label{lem-alg}
  Assume that $G$ is a Lie group with finitely many components.  Let
  $H\subseteq G_0$ be a connected closed normal subgroup of $G$ such
  that $H$ is locally algebraic, and such that $G/H$ is compact. Let
  $N$ denote the nilpotent radical of $H$, and let $S\subseteq N$ be a
  connected abelian normal subgroup of $G$. Let $\widehat{S}$ denote
  the character group of $S$ and let $G$ act on $\widehat{S}$ via
  conjugation.  Then the following assertions are true:
      \begin{enumerate}
      \item The orbit space $\widehat{S}/G$ is countably separated,
        i.e., all $G$-orbits in $\widehat{S}$ are locally closed.
      \item If $G_{\chi}$ is the stabilizer of some $\chi\in
        \widehat{S}$ for the action of $G$ on $\widehat{S}$, then
        $G_{\chi}/(G_{\chi})_0$ is amenable.
\end{enumerate}
\end{lem}
\begin{proof}
  We first show that it is sufficient to prove the result for the case
  $G=H$. Indeed, if we already know that $\widehat{S}/H$ is countably
  separated, then we observe that $\widehat{S}/H$ is a topological
  $G/H$-space such that $\widehat{S}/G\cong (\widehat{S}/H)/(G/H)$.
  But it is an easy exercise to prove that the quotient space of a
  countably separated space by a compact group action is countably
  separated.
  
  Assume now that $G_{\chi}$ is the stabilizer of some $\chi\in
  \widehat{S}$ in $G$. Then $H_{\chi}=G_{\chi}\cap H$ is the
  stabilizer in $H$. It follows that $H_{\chi}$ is a normal subgroup
  of $G_{\chi}$ such that $G_{\chi}/H_{\chi}$ is compact. If
  $H_{\chi}/(H_{\chi})_0$ is amenable, it also follows that
  $G_{\chi}/(H_{\chi})_0$, and hence also $G_{\chi}/(G_{\chi})_0$ are
  amenable.
  
  So, for the rest of the proof we assume that $G=H$.  In the next
  step we reduce to the case where $H$ is simply connected. For this
  let $\tilde{H}$ denote the universal covering group of $H$. Then
  $\tilde{H}$ has the same Lie algebra as $H$, and therefore it is
  locally algebraic. Let $q:\tilde{H}\to H$ denote the quotient map
  and let $C=\ker q$. Then $C$ is a discrete central subgroup of
  $\tilde{H}$.  Let $\frak{r}\subseteq \frak{h}$ denote the Lie
  algebra of $S$ and let $\tilde{S}$ denote the connected closed
  normal subgroup of $\tilde{H}$ corresponding to $\frak{r}$. Then
  $\tilde{S}$ is a vector subgroup of the nilpotent radical
  $\tilde{N}$ of $\tilde{H}$ and the quotient map $\tilde{H}\to H$
  maps $\tilde{S}$ surjectively onto $S$, i.e., we have $S\cong
  \tilde{S}/(\tilde{S}\cap C)$.  In particular, it follows that we may
  view $\widehat{S}$ as a closed $\tilde{H}$-invariant subspace of
  $\widehat{\tilde{S}}$, and we have
  $\widehat{S}/H=\widehat{S}/\tilde{H}$ (since the central subgroup
  $C$ acts trivially on $\widehat{S}$).  Thus, if
  $\widehat{\tilde{S}}/\tilde{H}$ is countably separated, the same is
  true for $\widehat{S}/H$.
  
  We now consider the stabilizers. It follows from the above
  considerations that if $H_{\chi}$ is the stabilizer of some $\chi\in
  \widehat{S}$, then $q^{-1}(H_{\chi})\subseteq \tilde{H}$ is the
  stabilizer of $\chi$ in $\tilde{H}$.  Thus it follows that
  $H_{\chi}=\tilde{H}_{\chi}/(C\cap \tilde{H}_{\chi})$.  Since the
  connected component of $\tilde{H}_{\chi}$ is mapped onto the
  connected component of $H_{\chi}$ under the quotient map, it follows
  that $H_{\chi}/(H_{\chi})_0$ is a quotient of
  $\tilde{H}_{\chi}/(\tilde{H}_{\chi})_0$.  Thus, if the latter is
  amenable, the same is true for $H_{\chi}/({H}_{\chi})_0$.
  
  Thus, in what follows we may assume without loss of generality that
  $H$ is simply connected. We are then in precisely the same situation
  as in the proof of Case (A) of the proof of the Theorem on page 2 of
  \cite{Puk}, and from now on we can follow the line of arguments as
  given on pages 2 and 3 of \cite{Puk} to see that $\widehat{S}/H$ is
  countably separated.  Moreover, the arguments presented in steps c)
  and d) on page 3 of Punk\'anzsky's book imply that for each
  stabilizer $H_{\chi}$ the quotient $H_{\chi}/(H_{\chi})_0$ is a
  finite extension of an abelian group, and hence is amenable.
\end{proof}

\begin{remark}\label{lem-stabilizers} We should point out that
  the result on the stabilizers in Lemma \ref{lem-alg} is most
  satisfying: Indeed if we know that every almost connected Lie group
  with dimension $\dim(G)< n$ satisfies BC for $\K$, say, then, by an
  easy application of Theorem \ref{thm-ext} the same is true for all
  Lie groups $H$ with $\dim(H)<n$ and $H/H_0$ amenable!
\end{remark}

Unfortunately, the result on the orbit space $\widehat{S}/G$ is not
sufficient for a direct application of Theorem \ref{thm-BCbundle}. So
we have to do some extra work to obtain more information on the
structure of $\widehat{S}/G$. To do this
we have to do two steps:
\begin{enumerate}
\item Reduce to cases where the action of $G$ on $\widehat{S}$
factors through an algebraic action of some real algebraic group
$G'$ (or a subgroup of finite index in $G'$).
\item Show that that the topological orbit-spaces of algebraic group actions 
on real affine varieties have nice stratifications
as required by Theorem \ref{thm-BCbundle}.
\end{enumerate}
Note that Puk\'anszky does the first reduction 
for the cocompact subgroup $H$ of $G$, which allowed us to
draw the conclusions of the previous lemma. 
However, with a bit more work we obtain similar conclusion 
for $G$. The following result is
certainly well-known to the experts, but since we didn't find a direct
reference we included the easy proof.

\begin{lem}\label{lem-red}
  Assume that $G$ is a Lie group with finitely many components such
  that $G$ has a connected closed normal subgroup $H$ with $H$
  semi-simple and $G_0/H$ a torus group.  Let $V$ be a finite
  dimensional real vector space and let $\rho:G\to \GL(V)$ be any
  continuous homomorphism. Then the Zariski closure $G'$ of $\rho(G)$
  is a (reductive) real algebraic group which contains $\rho(G)$ as a
  subgroup of finite index.
\end{lem}

\begin{proof}
 Let $R:\frak{g}\to\frak{gl}(V)$ denote
  the differential of $\rho$ and let $\frak{h}$ denote the ideal of
  $\frak{g}$ corresponding to $H$. Then $R(\frak{h})$ is semi-simple
  (or trivial). Since $\rho(H)$ is a semi-simple subgroup of $\GL(V)$
  it is closed in $\GL(V)$.  This follows from the fact that every
  semi-simple subalgebra of $\frak{gl}(V)$ is algebraic
(by Proposition \ref{prop-liealg}), which implies
  that $\rho(H)$ is the connected component of some algebraic linear
  subgroup of $\GL(V)$.  
Since $G/H$ is compact, it
  follows that $\rho(G)$ is a closed subgroup of $\GL(V)$, too.
  
  To simplify notation we assume from now on that $G$ itself is a
  closed subgroup of $\GL(V)$ and that $\rho$ is the
  identity map.  Let $\frak{g}=\frak{s}+\frak{z}$ be a Levi
  decomposition of $\frak{g}$. Since $G_0/H$ is abelian and $H$ is
  semi-simple, it follows that $\frak{s}=[\frak{g},\frak{g}]=\frak{h}$
  and $[\frak{h},\frak{z}]=[\frak{z},\frak{z}]=\{0\}$.  Thus
  $\frak{g}=\frak{h}\oplus\frak{z}$ is reductive.
  
  We now show that $\frak{g}$ is an algebraic subalgebra of $\frak{gl}(V)$.
By  Proposition \ref{prop-liealg} it suffices to show that 
$\frak{z}$ is algebraic and consists of semi-simple elements.
But this will follow if we can show that $Z=\exp(\frak{z})\subseteq 
\frak{gl}(V)$ is compact, and hence a torus group. 
Since $Z\cap H$ is finite (since every linear semi-simple group
  has finite center), the restriction to $Z$ of the quotient map
  $q:G\to G/H$ has finite kernel. Since $q(Z)=G_0/H$, $q(Z)$ is
  compact by assumption, and hence $Z$ is compact, too.
  
  It follows that the algebraic closure $\tilde{G}$ of $G_0$ is a
 reductive algebraic  group which contains $G_0$ as a subgroup of
  finite index. Since every element of $G$ fixes the Lie algebra
  $\frak{g}$ via the adjoint action, it also normalizes $\tilde{G}$.
  Therefore, $G'=G\tilde{G}$ is a reductive algebraic group which
  contains $G$ as a subgroup of finite index.
\end{proof}

We now show that quotients of linear algebraic group actions 
on affine varieties have nice stratifications in the sense of 
Theorem \ref{thm-BCbundle}. We are very grateful to J\"org Sch\"urmann
and Peter Slodowy for some valuable comments, which helped us 
to replace a  previous version of the following result
(which, as was pointed out to us by J\"org Sch\"urmann, contained a gap)
by

\begin{prop}\label{prop-reductive}
    Suppose that $G$ is a closed subgroup of finite index
 of a Zariski closed subgroup $G'$ of 
    $\GL(n,\RR)$ and that $V\subseteq \RR^n$ is a $G'$-invariant Zariski 
    closed subset of $\RR^n$. Then there exists a 
    stratification 
    $$\emptyset=V_0\subseteq V_1\subseteq\cdots \subseteq V_l=V$$
    of open $G$-invariant subsets $V_i$ of $V$
    such that 
    $(V_i\setminus V_{i-1})/G$ admits a continuous and open finite-to-one map
    onto a differentiable manifold.
\end{prop}

Since every manifold has a triangulation, the above result
really gives what we need to apply Theorem \ref{thm-BCbundle}.
For the proof we need the following lemma about certain decompositions of 
continuous semi-algebraic maps.

\begin{lem}\label{lem-semialg}
    Let $X,Y$ be semi-algebraic sets and let $f:X\to Y$ 
    be a continuous semi-algebraic map (see \cite{BCR} for 
    the notations).
    Then there exists a stratification
    $$\emptyset= Z_0\subseteq Z_1\cdots \subseteq Z_l=f(X),$$
    with each $Z_i$ open in $f(X)$, $Z_i\setminus Z_{i-1}$ 
    is a differentiable manifold and 
    $$f:f^{-1}(Z_i\setminus Z_{i-1})\to Z_i\setminus Z_{i-1}$$
    is open (in the euclidean topology) for all $1\leq i\leq l$.
\end{lem}
\begin{proof} Since the image of a semi-algebraic set under a 
    semi-algebraic map  is 
    semi-algebraic (see \cite[Proposition 2.2.7]{BCR}),
    we may assume without loss of generality that $Y=f(X)$.
    By \cite[Corollary 9.3.3]{BCR} there exists a closed 
    semi-algebraic subset $Y_1\subseteq Y$ with $\dim(Y_1)<\dim(Y)$,
    such that $Y\setminus Y_1$ is a finite disjoint
    union of connected components (combine with \cite[Theorem 
    2.4.5]{BCR}) and such that the restriction of $f$ to the inverse 
    image of each component is a projection, hence open.
    Thus $f:f^{-1}(Y\setminus Y_1)\to Y\setminus Y_1$ is open, too.
Indeed, the construction (using \cite[Proposition 9.18]{BCR})
implies that $Y\setminus Y_1$ is homeomorphic to 
a submanifold of some $\RR^m$.
Put $Z_0=Y\setminus Y_1$. Since $\dim(Y_1)<\dim(Y)$, the result
    follows by induction.
\end{proof}

\begin{remark}\label{rem-orbits} 
    Let $G\subseteq \GL(n,\RR)$ be a real linear algebraic group,
    and let $G_{\CC}\subseteq \GL(n,\CC)$ be its complexification.
    Then it follows from \cite[Proposition 2.3]{BHC} that
     each $G_{\CC}$-orbit in $\CC^n$ contains at most finitely many
    $G$-orbits in $\RR^n\subseteq\CC^n$. 
\end{remark}
    
\begin{proof}[Proof of Proposition \ref{prop-reductive}]
We first note that we may assume without loss of generality 
that $G=G'$. Indeed, since $G$ has finite index in $G'$,
every $G'$-orbit decomposes into finitely many $G$-orbits.
Thus, if $\emptyset=V_0\subseteq V_1\subseteq\cdots \subseteq V_l=V$
is a stratification of $V$ for the $G'$-action
with the required properties, 
it is also a stratification for the $G$-action with the same properties.
Thus we assume from now on that $G$ is a Zariski closed subgroup of 
$\GL(V)$.

    Let $V_{\CC}\subseteq \CC^n$ denote the 
    complexification 
    of $V$. Consider the diagram
    $$
    \begin{CD}
        V @>>> V_{\CC}\\
        @VVV      @VVV\\
        V/G  @>>>  V_{\CC}/G_{\CC}.
\end{CD}
$$
 By the theorem of Rosenlicht (\cite{R}, but see also 
    \cite[Satz 2.2 on p. 23]{Kr}), there exists a sequence
    $$V_{\CC}=W_0\supseteq W_1\supseteq W_2\supseteq \cdots \supseteq 
    W_r=\emptyset,$$
    of Zariski-closed $G_{\CC}$-invariant subsets of 
    strictly decreasing dimension such that $W_i\setminus W_{i+1}$
    has closed $G_{\CC}$-orbits and the geometric quotient by 
    $G_{\CC}$
    of $W_i\setminus W_{i+1}$ exists. This means that the
    quotient $(W_i\setminus W_{i+1})/G_{\CC}$ 
    can be realized as an algebraic set and the quotient map is
    also algebraic.
    Let $\mathcal O$ be the first of the sets $W_i\setminus W_{i+1}$ 
    which has nonempty intersection with $V$. Restricting the maps 
    in the above diagram gives 
     $$
    \begin{CD}
        V\cap\mathcal O @>>> \mathcal O\\
        @VVV      @VVV\\
        (V\cap \mathcal O)/G  @>>>  \mathcal O/G_{\CC}.
\end{CD}
$$
The resulting map $f$ from $V\cap\mathcal O$ to $\mathcal O/G_{\CC}$ is 
an algebraic map, and hence it is a continuous semi-algebraic 
map. Thus it follows from Lemma \ref{lem-semialg} that, if $Y$ denotes the 
image of $X:=V\cap \mathcal O$ in $\mathcal O/G_{\CC}$, then 
$Y$ has a stratification 
$$\emptyset = Z_0\subseteq Z_1\subseteq \cdots \subseteq Z_s= Y$$
such that $f:f^{-1}(Z_i\setminus Z_{i-1})\to Z_i\setminus Z_{i-1}$
    is open for all $1\leq i\leq s$, each $Z_i$ is open in 
    $Y$, and the difference sets $Z_i\setminus Z_{i-1}$ are 
submanifolds of some $\RR^m$.
Put $V_i=f^{-1}(Z_i)$ for $0\leq i\leq s$. Then $V_s=V\cap \mathcal O$.
By Remark \ref{rem-orbits}, if we pass through the lower left corner of the 
diagram, the corresponding maps  $(V_i\setminus V_{i-1})/G\to 
Z_i\setminus Z_{i-1}$ are open, finite-to-one, onto the 
manifolds $Z_i\setminus Z_{i-1}$.

Now replace $V$ by the invariant Zariski-closed subset 
$V\setminus\mathcal O$. Repeating the above arguments finitely many
times gives the desired stratification (the procedure stops 
after finitely many steps, since any increasing sequence of 
Zariski open sets eventually stabilizes).
\end{proof}

Using the above results, we are now able to prove

\begin{prop}\label{prop-reduce}
  Suppose that $G$ is a Lie group with finitely many components and with
  connected closed normal subgroups $N\subseteq H\subseteq
  G_0\subseteq G$ as in (\ref{eq-structure}), i.e., $N$ is the
  nilradical of $H$, $H/N$ is semi-simple and $G_0/H$ is a torus group.
  Let $S\subseteq Z(N)$ be a connected closed subgroup which is normal
  in $G$, where $Z(N)$ denotes the center of $N$.  Then 
$\widehat{S}$ decomposes into a countable disjoint union of 
open $G$-invariant sets $V_n$ such that each $V_n$ has a
  stratification
  $$\emptyset=U_{0}\subseteq U_{1}\subseteq\cdots \subseteq
  U_l=V_n$$
(where $l$ may depend on $n$)
  of open $G$-invariant subsets of $V_n$,
  and continuous open surjections
  $$q_i:U_i\smallsetminus U_{i-1}\to Y_i,\quad 1\leq i\leq l,$$
  such
  that each $Y_i$ is a differentiable manifold and inverse images of
  points in $Y_i$ are finite unions of $G$-orbits in $V_n$ for
  all $1\leq i\leq l$.
\end{prop}
\begin{proof}
  Let $\frak{s}$ denote the ideal of $\frak{g}$ corresponding to $S$.
  Then we may identify $\widehat{S}$ with a closed ${G}$-stable subset
  of $\frak{s}^*$ of the form $R\times Z$ with $R$ being a vector
  subgroup of $\frak{s}^*$ and $Z$ a finitely generated free abelian
  group.  Note that $Z$ can be identified with the dual of the maximal
  compact subgroup in $S$, and therefore we can decompose $Z$ into a
  disjoint union of $G$-orbits, which are all finite since $G_0$ acts
  trivially on $Z$.  It then follows that $\widehat{S}$ can be
  decomposed into a disjoint union of $G$-invariant sets of the form
  $R\times F$ with $F\subseteq Z$ finite.
  
  The action of ${G}$ on $\widehat{S}$ is given via the coadjoint
  representation $\Ad^*_\frak{s}:{G}\to GL(\frak{s}^*)$.  Since
  $S\subseteq Z({N})$, it follows that this representation factors
  through a representation of $G/N$.  Thus it follows from the general
  assumptions on $G$ and Lemma \ref{lem-red} that the algebraic
  closure $G'$ of $\Ad^*_\frak{s}(G)$ in $GL(\frak{s}^*)$ is a
  reductive algebraic group which contains the image of
  $\Ad^*_\frak{s}({G})$ as a subgroup of finite index.  Since the
  $G$-stable sets of the form $R\times F$ of the previous paragraph
  are closed algebraic subvarieties of $\frak{s}^*$, it follows that
  these sets are also invariant under the action of the Zariski
  closure $G'$ of $\Ad^*_\frak{s}(G)$. Thus it follows from
  Proposition \ref{prop-reductive} that for each such set we obtain a
  stratification
  $$\emptyset =U_0\subseteq U_1\subseteq \cdots \subseteq
  U_l=R\times F$$
with the required properties.
\end{proof}

We are now ready for the final step:

\begin{proof}[Proof of Proposition \ref{prop-main}]
  By the discussion at the beginning of this section we may assume
  without loss of generality that $G$ is as in (\ref{eq-structure}),
  i.e., we have connected closed normal subgroups
  $$N\subseteq H\subseteq G_0\subseteq G$$
  such that $N$ is a
  nontrivial nilpotent group $H/N$ is semi-simple and $G_0/H$ is a
  torus group.  For the induction step we have to show that every
  central extension
  $$1\to \TT\to\bar{G}\to G\to 1$$
  satisfies BC for $\CC$.  Let
  $\bar{N}$, $\bar{H}$ and $\bar{G}_0$ denote the inverse images of
  $N$, $H$ and $G_0$ in $\bar{G}$.  Then the sequence of normal
  subgroups
  $$\bar{N}\subseteq \bar{H}\subseteq \bar{G}_0\subseteq\bar{G}$$
  has
  the same general properties as the sequence $N\subseteq H\subseteq
  G_0\subseteq G$, in particular, $\bar{N}$ is the nilradical of
  $\bar{H}$ and $\bar{H}$ is locally algebraic. Let $T$ denote the
  central copy of $\TT$ in $\bar{G}$ coming from the given central
  extension. We now divide the proof into the following cases:
\begin{enumerate}
\item[C(1)] The center $S=Z(\bar{N})$ of $\bar{N}$ has dimension
  greater or equal to two.
\item[C(2)] $Z(\bar{N})=T$.
\end{enumerate}

We start with Case C(1): By Theorem \ref{thm-ext} (and the discussion
following that theorem) it suffices to show that $\bar{G}/S$ satisfies
BC with coefficients in $C_0(\widehat{S},\K)$ where the action of
$\bar{G}/S$ on $\widehat{S}$ is given by conjugation.  By Theorem
\ref{thm-BCbundle} it suffices to show that all stabilizers
$(\bar{G}/S)_{\chi}=\bar{G}_{\chi}/S$ satisfy BC for $\K$ and that
$\widehat{S}$ has a nice stratification.  While the latter follows
from Proposition \ref{prop-reduce}, the requirement on the stabilizers
follows from Lemma \ref{lem-alg}, Remark \ref{lem-stabilizers}, and the
induction assumption since
$$\dim(\bar{G}_{\chi}/S)\leq \dim(\bar{G})-2<\dim(G).$$
This finishes
the proof in Case C(1).

For the proof of Case C(2) we have to do some more reduction steps in
order to use the same line of arguments as in C(1).  For this it is
useful to consider the following two subcases:

\begin{enumerate}
\item[(2)a] If $\bar{Z}(N)$ denotes the inverse image of the center
  $Z(N)$ of $N$ in $\bar{G}$, then $Z(\bar{Z}(N))=T$.
\item[(2)b] $\dim\big(Z(\bar{Z}(N))\big)\geq 2$.
\end{enumerate}

In Case (2)a we consider the normal subgroup $S=\bar{Z}(N)$ of $G$.
Then $S$ is a connected two-step nilpotent Lie group with
one-dimensional center $T$, and therefore a Heisenberg group.  It
follows that $C^*_r(S)=C^*(S)$ can be written as the direct sum
$$C^*(S)=\oplus_{\chi\in \widehat{T}}A_{\chi}$$
with
$$A_{\chi}\cong \K,\quad\text{if $\chi\neq 1$, and}\quad
A_1=C_0(\widehat{S/T}).$$
Since $\bar{G}$ acts trivially on
$\widehat{T}$, it follows that the decomposition action of $\bar{G}/S$
on $C^*(S)\otimes \K$ induces an action on each fibre $A_{\chi}$, and,
by Theorem \ref{thm-ext} together with Proposition \ref{prop-limit},
it follows that $\bar{G}$ satisfies BC with coefficients in $\CC$ if
$\bar{G}/S$ satisfies BC with coefficients in $A_{\chi}\otimes\K$ for
each $\chi\in \widehat{T}$.  If $\chi\neq 1$, we get
$A_{\chi}\otimes\K\cong \K$, and the desired result follows from the
induction assumption and the fact that $\dim(\bar{G}/S)<\dim(G)$.

So we only have to deal with the case $\chi=1$, where we have to deal
with the fibre $C_0(\widehat{S/T},\K)=C_0(\widehat{Z(N)},\K)$.  But
here we are exactly in the same situation as in the proof of Case C(1),
since the action of $\bar{G}$ on $\widehat{Z(N)}$ factors through an
action of $G/N$ and all stabilizers of the characters of $Z(N)$ have
dimension strictly smaller than $\dim(G)$.

\medskip

We have to work a bit more for the Proof of Case (2)b.  Here we put
$S=Z(\bar{Z}(N))$. Then $S$ is a connected abelian subgroup of
$\bar{N}$ and it follows from Lemma \ref{lem-alg}, Remark
\ref{lem-stabilizers}, the fact that $\dim(\bar{G}/S)<\dim(G)$ 
and the induction assumption that
all stabilizers for the action of $\bar{G}/S$ on $\widehat{S}$ satisfy
BC for $\K$.

Again we study the structure of the orbit space $\widehat{S}/\bar{G}$.
For each $\chi\in \widehat{T}$ we define
$$\widehat{S}_{\chi}=\{\mu\in \widehat{S}: \mu|_T=\chi\}.$$
Since $T$
is central in $\bar{G}$, it follows that $\bar{G}$ acts trivially on
$\widehat{T}$, and hence that $\widehat{S}_\chi$ is
$\bar{G}$-invariant for all $\chi\in \widehat{T}$.  Since
$\widehat{T}$ is discrete, we may write
$$C_0(\widehat{S},\K)\cong \oplus_{\chi\in
  \widehat{T}}C_0(\widehat{S}_{\chi},\K)$$
with fiberwise action of
$\bar{G}/S$.  Thus by continuity of BC it suffices to deal with the
single fibers. For $\chi=1$ we are looking at the action of
$\bar{G}/S\cong G/(S/T)$ on $\widehat{S}_1\cong \widehat{S/T}$, and
since $S/T$ is a central subgroup of $N$ we may again argue precisely
as in the proof of Case C(1) to see that $\bar{G}/S$ satisfies BC for
$C_0(\widehat{S}_1,\K)$.

In order to deal with the other fibers we are now going to show that
$\bar{G}$ acts transitively on $\widehat{S}_\chi$ for each nontrivial
character $\chi\in \widehat{T}$.  It follows then directly from
Corollary \ref{cor-induced} that $\bar{G}/S$ satisfies BC for
$C_0(\widehat{S}_\chi,\K)$.  In fact, Lemma \ref{lem-nilpotent} below
shows that $\bar{N}$ already acts transitively on $\widehat{S}_{\chi}$
for $\chi\neq 1$ and the result will follow from that lemma.
\end{proof}

The following lemma is certainly well known to the experts on the
representation theory of nilpotent groups. For the readers convenience
we give the elementary proof.

\begin{lem}\label{lem-nilpotent}
  Assume that $N$ is a connected nilpotent Lie group with
  one-dimensional center $Z(N)=T$. Let $S$ be a closed connected
  abelian normal subgroup of $N$ such that $T\subseteq S$ and
  $S/T\subseteq Z(N/T)$. Let $1\neq \chi\in \widehat{T}$ and let
  $\widehat{S}_{\chi}=\{\mu\in \widehat{S}:\mu|_T=\chi\}$.  Then $N$
  acts transitively on $\widehat{S}_\chi$ by conjugation.
\end{lem}
\begin{proof} We may assume without loss of generality that $N$ is simply
  connected. In fact, if this is not the case, we pass to the
  universal covering group $\tilde{N}$ of $N$ and the universal
  covering $\tilde{S}\subseteq \tilde{N}$ of $S$ and observe that
  there exists a discrete subgroup $D\subseteq \tilde{T}=Z(\tilde{N})$
  such that $N=\tilde{N}/D$, $S=\tilde{S}/D$, $T=\tilde{T}/D$ and
  $\widehat{S}_{\chi}$ can then be (equivariantly) identified with
  $\widehat{\tilde{S}}_{\chi}$ for all $\chi\in
  \widehat{T}\subseteq\widehat{\tilde{T}}$.
  
  Let $\frak{n}$, $\frak{s}$ and $\frak{t}$ denote the Lie algebras of
  $N$, $S$ and $T$, respectively. Since $N$ is simply connected, we
  can write
  $$N=\{\exp(X): X\in \frak{n}\}$$
  with multiplication given by the
  Campbell-Hausdorff formula.  In particular, if $Y\in \frak{s}$, then
  $$\exp(X)\exp(Y)=\exp(X+Y+[X,Y])$$
  for all $X\in \frak{n}$, since it
  follows from the assumption that $S/T\subseteq Z(N/T)$ that
  $[X,Y]\in \frak{t}=\frak{z}(\frak{n})$ and all commutators with
  $[X,Y]$ vanish. In particular, if we conjugate $\exp(Y)$ by
  $\exp(X)$ we get the formula
\begin{equation}\label{eq-conj}
\exp(X)\exp(Y)\exp(-X)= \exp(Y+[X,Y])
\end{equation}
for all $Y\in \frak{s}$.

Assume now that $\dim(\frak{s})=n+1$ and let $0\neq Z\in \frak{t}$.
We claim that we can find a basis $\{Y_1,\ldots,Y_n, Z\}$ of
$\frak{s}$ and elements $X_1,\ldots, X_n\in \frak{n}$ such that
\begin{equation}\label{eq-schmidt}
      [X_i,Y_i]=Z\quad\text{and}\quad [X_i, Y_j]=0
\end{equation}
for all $1\leq i,j\leq n, i\neq j$.  Indeed this follows from an easy
Schmidt-orthogonalization procedure: Consider the bilinear form
$$(\cdot, \cdot):\frak{n}\times \frak{s}\to \RR;\quad (X,Y)=\lambda
\Leftrightarrow [X,Y]=\lambda Z.$$
Choose $Y_1\in
\frak{s}\smallsetminus\frak{t}$. Then there exists $X_1\in \frak{n}$
with $[X_1,Y_1]\neq 0$ (since $Y_1$ is not central), and hence we may
assume that $[X_1,Y_1]=Z$. Assume now that we have chosen $Y_1,
\ldots, Y_l$ and $X_1,\ldots , X_l$, $l<n$, with the desired
properties. Choose $Y_{l+1}$ not in the span of $\{Y_1,\ldots,
Y_l,Z\}$ and $X'\in \frak{n}$ with $[X', Y_{l+1}]=Z$. Then it is easy
to check that $X_{l+1}=X'-\sum_{i=1}^l(X',Y_i)Y_i$ satisfies
$[X_{l+1},Y_{l+1}]=Z$ and $[X_{l+1}, Y_i]=0$ for $1\leq i\leq l$.

Now we identify $S$ with $\frak{s}$ (via $\exp$) and $\widehat{S}$
with $\frak{s}^*$. The conjugation action of $N$ on $\widehat{S}$ is
then transferred to the coadjoint action $\Ad^*$. If $\{f_1,\ldots,
f_n, g\}$ is a dual basis for the basis $\{Y_1,\ldots, Y_n, Z\}$ of
$\frak{s}$, the result will follow if we can show that
$$\Ad^*(N)(\lambda g)=\sp\{f_1,\ldots, f_n\}+\lambda g$$
for all
$0\neq\lambda\in \RR$. By rescaling we may assume that $\lambda=1$.
But for $\lambda_1,\ldots,\lambda_n\in \RR$ we can compute
$$\Big(\Ad^*(\exp(\lambda_1X_1+\cdots+\lambda_n X_n))(g)\Big)(Y_i)
=g(Y_i+\lambda_i[X_i,Y_i])=\lambda_ig(Z)=\lambda_i.$$
Since $Z$ is
central in $\frak{n}$ it follows that
$\Big(\Ad^*(\exp(X))(g)\Big)(Z)=g(Z)=1$ for all $X\in\frak{n}$.  Thus
$$\Ad^*(\exp(\lambda_1X_1+\cdots+\lambda_n X_n))(g)=
\lambda_1f_1+\cdots +\lambda_nf_n+g.$$
\end{proof}

\section{Relations to the $K$-theory of the maximal compact subgroup}
\label{sec-max}

In this section we want to describe the relations between the
$K$-theory of $C_r^*(G)$ and the $K$-theory of $C^*(L)$, where $L$
denotes the maximal compact subgroup of the almost connected group
$G$ (we chose the letter $L$ to avoid confusion).  We should mention
that all results presented here (accept the conclusions drawn out of
our main theorem) are well known, but since they have important impact
on our results, we found it useful to give at least a brief report.
The main references for these results are \cite{CM, Kas0}, and we
refer especially to \cite[\S4]{CM} for a more geometric discussion of
some the results presented in this section.

If $G$ and $L$ are as above, it follows from work of Abels (see
\cite{Ab}) that $G/L$ is a universal proper $G$-space.  Thus we have
$$\k_*^{\top}(G,A)\cong \KK_*^G(C_0(G/L),A)\stackrel{\res^G_L}{\cong}
\KK_*^L(C_0(G/L),A),$$
where the second isomorphism follows from
\cite[Corollary to Theorem 5.7]{Kas2}.  Also by the work of Abels
\cite{Ab}, $G/L$ is a Riemannian manifold which is $L$-equivariantly
diffeomorphic to the tangent space $V:=T_{eL}$ equipped with the
adjoint action of $L$ on $V$. It follows then from Kasparov's work in
\cite{Kas0} (see \cite[Lemma 7.7]{CE1} for a more extensive
discussion) that tensoring with $C_0(V)$ gives a natural isomorphism
$$\sigma_{C_0(V)}:\KK_*^L(C_0(V),A)\to \KK_*^L(C_0(V)\otimes C_0(V),
A\otimes C_0(V)),$$
and by Kasparov's Bott-periodicity theorem (see
\cite[Theorem 7]{Kas0}) we know that $C_0(V)\otimes C_0(V)$, equipped
with the diagonal action, is $\KK^L$-equivalent to $\CC$ (but see also
the discussion below).  Thus we obtain the following chain of
isomorphisms
\begin{align*}
  \KK_*^L(C_0(G/L),A)&\cong
  \KK_*^L(C_0(V),A)\stackrel{\sigma_{C_0(V)}}{\cong}
  \KK_*^L(C_0(V){\otimes}C_0(V), A{\otimes}C_0(V))\\
  &\cong \KK_*^L(\CC, A{\otimes}
  C_0(V))=\k_*\big((A{\otimes}C_0(V))\rtimes L\big),
\end{align*}
where the last isomorphism follows from the Green-Julg theorem.
Hence, as a direct consequence of Theorem \ref{thm-main1} we can
deduce

\begin{thm}\label{thm-connected}
  Assume that $G$ is an almost connected (second countable) group with
  maximal compact subgroup $L$.  Let $\K=\K(H)$ be the algebra of
  compact operators on the separable Hilbert space $H$ equipped with
  any action of $G$.  Then $\k_*(\K\rtimes_rG)$ is naturally
  isomorphic to $\k_*\big((\K{\otimes}C_0(V))\rtimes L\big)$.
\end{thm}

By Kasparov's Bott-periodicity theorem (see \cite[Theorem 7]{Kas0}) it
follows that $C_0(V)$ is $\KK^L$-equivalent to the graded complex
Clifford algebra $Cl(V)$ (with respect to a compatible inner product
on $V$), equipped with the action of $L$ induced by the given action
on $V$.  So we can replace $C_0(V)$ by the graded $C^*$-algebra
$Cl(V)$, but then we have to use graded $K$-theory!

Let us look a bit closer to the implications of this Bott-periodicity
theorem. Assume for the moment that $V$ is even dimensional and that
the action of $L$ on $V$ preserves a given orientation of $V$, i.e.,
the action factors through a homomorphism $\varphi:L\to SO(V)$.  We
have a central extension
$$0\to \TT\to \Spin^c(V)\to SO(V)\to 0$$
of $SO(V)$, where
$\Spin^c(V)\subseteq CL(V)$ denotes the group of complex spinors (e.g.
see \cite{ABS}).  The corresponding action of $L$ on $Cl(V)$ is given
by the homomorphism
$$L\to SO(V)\cong \Spin^c(V)/\TT=\Ad(\Spin^c(V)).$$
Now choose a fixed
orthonormal base $\{e_1,\ldots,e_n\}$ of $V$.  Then the grading of
$Cl(V)$ is given by conjugation with the symmetry $J=e_1\cdot
e_2\cdots e_n\in Cl(V)$. One can show that, up to a sign, $J$ does not
depend on the choice of this basis, and the sign only depends on the
orientation of the basis. In particular, $J$ is invariant under
conjugation with elements in $\Spin^c(V)$.  From this it follows that
the graded $L$-algebra $Cl(V)$ is $L$-equivariantly Morita equivalent
to the trivially graded $L$-algebra $Cl(V)$ -- a Morita equivalence is
given by the module $Cl(V)$ with given $L$-action and grading
automorphism given by left multiplication with $J$.  Moreover, since
$n=\dim(V)$ is even, $Cl(V)$ is isomorphic to the simple matrix
algebra $M_{2^n}(\CC)$.

Assume now that $\dim(G/L)$ is odd. Then, replacing $G$ by
$G\times\RR$ (with trivial action of $\RR$ on $\K$) we get
$$\k_*(\K\rtimes_rG)=\k_{*+1}(\K\rtimes_r(G\times\RR)).$$
Moreover, if
the action of $L$ on $V=T_{eL}$ is orientation preserving, the same is
true for the resulting action of $L$ on $V\times \RR$, which we
identify with the tangent space at $eL$ in the group $G\times\RR$.
Hence, modulo a dimension shift, we can use the above considerations
also for this case. Thus, as a consequence of Theorem
\ref{thm-connected} we obtain

\begin{thm}\label{thm-cliff}
  Assume that $G$ is an almost connected group with maximal compact
  subgroup $L$ such that the adjoint action of $L$ on $V=T_{eL}$ is
  orientation preserving. Then there are natural isomorphisms
  $$\k_*(\K\rtimes_rG)\cong \k_*\big((\K\otimes Cl(V))\rtimes L\big)$$
  if $\dim(G/L)$ is even and
  $$\k_{*+1}(\K\rtimes_rG)\cong\k_*\big((\K\otimes
  Cl(V\times\RR))\rtimes L\big)$$
  if $\dim(G/L)$ is odd. Here all
  algebras are trivially graded!
\end{thm}

Perhaps, the above result has its most satisfying formulation if
translated into the language of twisted group algebras.  For this let
$\om\in Z^2(G,\TT)$ denote a representative of the Mackey obstruction
for the action of $G$ on $\K$ (see the discussion preceding Lemma
\ref{lem-twist}).  Then $\K\rtimes_rG$ is isomorphic to
$C_r^*(G,\om)\otimes\K$, where $C_r^*(G,\om)$ denotes the reduced
twisted group algebra $C_r^*(G,\om)$ (e.g., see \cite[Theorem
18]{green1}).  Recall that $C^*_r(G,\om)$ can be defined either as the
reduced twisted crossed product $\CC\rtimes_r(G_{\om},\TT)$ with
respect to the twisted action $(\id,\chi_1)$ (which, by Lemma
\ref{lem-twist}, is Morita equivalent to the given action on $\K$), or
as the completion of $L^1(G)\subseteq B(L^2(G))$, where $L^1(G)$ acts
on $L^2(G)$ by the twisted convolution
$$f*\xi(s)=\int_Gf(t)\om(t,t^{-1}s)\xi(t^{-1}s)\,dt,\quad f\in L^1(G),
\xi\in L^2(G).$$
Up to isomorphism, $C_r^*(G,\om)$ only depends on the
class $[\om]\in H^2(G,\TT)$.  Conversely, given any cocycle, the
representation $\lambda_{\om}:G\to U(L^2(G))$ given by
$$\big(\lambda_{\om}(t)\xi\big)(s)=\om(t,t^{-1}s)\xi(t^{-1}s)$$
determines an action of $G$ on $\K(L^2(G))$ with Mackey obstruction
represented by $\om$.

Note that the Mackey obstruction for the action of $L$ on $\K$ is
given by the restriction of $\om$ to $L$ and the obstruction for the
action of $L$ on $Cl(V)\cong M_{2^n}(\CC)$ (if $\dim(V)$ is even) is
given by the pull-back, say $\mu_L$, to $L$ of a cocycle representing
the central extension
$$1\to \TT\to\Spin^c(V)\to SO(V)\to 1.$$
Since $\Spin^c(V)\cong
(\TT\times \Spin(V))/\ZZ_2$ (diagonal action), where
$$1\to \ZZ_2\to \Spin(V)\to SO(V)\to 1$$
is the real group of spinors,
the cocycle $\mu_L$ can be chosen to take values in the subgroup
$\ZZ_2\subseteq \TT$, and therefore $\mu_L^2=1$.  Note that $\mu_L$ is
trivial if and only if the homomorphism $\varphi:L\to SO(V)$
factorizes through $\Spin^c(V)$ (i.e., if and only if $G/L$ carries a
$G$-invariant $\Spin^c$-structure).  If $\dim(V)$ is odd, we may
define $\mu_L$ in the same way as above, noticing that this cocycle is
equivalent to the pull back of (a cocycle representing) the extension
$$1\to\TT\to\Spin^c(V\times\RR)\to SO(V\times\RR)\to 1,$$
which
follows from the fact that $L$ acts trivially on $\RR$!  Since the
Mackey obstruction of a tensor product of actions is the product of
the Mackey obstructions of the factors, we obtain

\begin{thm}\label{thm-twist}
  Assume that $G$ is an almost connected group with maximal compact
  subgroup $L$ such that the adjoint action of $L$ on $V=T_{eL}$ is
  orientation preserving. Let $n=\dim(G/L)$ and let $\om\in
  Z^2(G,\TT)$ be any cocycle on $G$.  Then
  $$\k_*\big(C_r^*(G,\om)\big)\cong \k_{*+n}\big(C^*(L,
  \om\cdot\mu_L)\big).$$
  In particular, in the special case where
  $\om$ is trivial, we obtain an isomorphism
  $$\k_*\big(C_r^*(G)\big)\cong \k_{*+n}\big(C^*(L, \mu_L)\big).$$
\end{thm}

Again, $\mu_L$ is trivial if and only if $G/L$ carries a $G$-invariant
$\Spin^c$-structure.  In general, since $C^*(L,\om\cdot\mu_L)$ is the
quotient of the central extension $L_{\om\cdot\mu_L}$ of $L$ by $\TT$
corresponding to the character $\chi_1$ of $\TT$, it follows that
$C^*(L,\om\cdot\mu_L)$ is a direct sum of (possibly infinitely many)
matrix algebras.  Thus as a direct corollary of the above result we
obtain:

\begin{cor}\label{cor-kgroup}
  Assume that $G$, $L$ and $\om$ are as in Theorem \ref{thm-twist}.
  Then $\k_{0+n}\big(C_r^*(G,\om)\big)$ is isomorphic to a free
  abelian group in at most countably many generators and
  $\k_{1+n}\big(C_r^*(G,\om)\big)=\{0\}$.
\end{cor}

This result has interesting consequence towards the question of
existence of square integrable representations of connected unimodular
Lie groups. In fact, combining the above corollary with \cite[Theorem
4.6]{ros} gives:

\begin{cor}[{cf \cite[Corollary 4.7]{ros}}]\label{cor-square}
  Let $G$ be a connected unimodular Lie group.  Then all
  square-integrable factor representations of $G$ are type I.
  Moreover, $G$ has no square-integrable factor representations if
  $\dim(G/L)$ is odd, where $L$ denotes the maximal compact subgroup
  of $G$.
\end{cor}

We refer to \cite{ros} for more detailed discussions on this kind of
applications of the positive solution of the Connes-Kasparov
conjecture.  Note that Theorem \ref{thm-twist} and Corollary
\ref{cor-kgroup} do not hold in general without the assumption that
the action of $L$ on $V=T_{eL}$ is orientation preserving. In fact an
easy six-term-sequence argument shows that it cannot hold for the
group $G=\RR\rtimes \ZZ_2$, where $\ZZ_2$ acts on $\RR$ by reflection
through $0$.

\def\mathcs{{\normalshape\text{C}}^{\displaystyle *}}

\end{document}